\documentclass{article}

\usepackage{arxiv}

\usepackage[utf8]{inputenc} 
\usepackage[T1]{fontenc}    
\usepackage{hyperref}       
\usepackage{url}            
\usepackage{booktabs}       
\usepackage{amsfonts}       
\usepackage{nicefrac}       
\usepackage{microtype}      
\usepackage{lipsum}
\usepackage{lscape}
\usepackage[authoryear, sort]{natbib}

\usepackage{morefloats}
\usepackage{color}
\usepackage{multirow}
\usepackage{latexsym}
\usepackage{amsmath,amssymb,amsthm,graphicx}
\usepackage{dsfont}
\usepackage{enumitem}

\newtheoremstyle{thm}
{9pt}
{9pt}
{\itshape}
{}
{\bfseries}
{.}
{ }
{}
\theoremstyle{thm}

\newtheorem{theorem}{Theorem}[section]
\newtheorem{lemma}[theorem]{Lemma}
\newtheorem{corollary}[theorem]{Corollary}

\newtheoremstyle{def}
{9pt}
{9pt}
{}
{}
{\bfseries}
{.}
{ }
{}
\theoremstyle{def}

\newtheorem{remark}[theorem]{Remark}

\newenvironment{prf}{\textbf{\emph{Proof.}}}{\qed}

\newcommand{\R}{\mathbb{R}} 
\newcommand{\N}{\mathbb{N}} 
\newcommand{\E}{\mathbb{E}} 

\renewcommand{\footnoterule}{%
	\kern -3.5pt
	\hrule width \textwidth height 1pt
	\kern 3.5pt
}

\makeatletter
\def\blfootnote{\xdef\@thefnmark{}\@footnotetext}
\makeatother

\title{Minimum $L^q$-distance estimators for non-normalized parametric models}

\author{S. Betsch\\
Institute of Stochastics,\\ Karlsruhe Institute of Technology (KIT),\\ Germany.\\
\href{mailto:Steffen.Betsch@kit.edu}{Steffen.Betsch@kit.edu}\\
\And  B. Ebner\\
Institute of Stochastics,\\ Karlsruhe Institute of Technology (KIT),\\ Germany.\\
\href{mailto:Bruno.Ebner@kit.edu}{Bruno.Ebner@kit.edu}\\
\And B. Klar\\
Institute of Stochastics,\\ Karlsruhe Institute of Technology (KIT),\\ Germany.\\
\href{mailto:Bernhard.Klar@kit.edu}{Bernhard.Klar@kit.edu}\\
}

\begin{document}

\date{\today}
\maketitle


\begin{abstract}
 We propose and investigate a new estimation method for the parameters of models consisting of smooth density functions on the positive half axis. The procedure is based on a recently introduced characterization result for the respective probability distributions, and is to be classified as a minimum distance estimator, incorporating as a distance function the $L^q$-norm. Throughout, we deal rigorously with issues of existence and measurability of these implicitly defined estimators. Moreover, we provide consistency results in a common asymptotic setting, and compare our new method with classical estimators for the exponential-, the Rayleigh-, and the Burr Type XII distribution in Monte Carlo simulation studies. We also assess the performance of different estimators for non-normalized models in the context of an exponential-polynomial family.
\end{abstract}

\section{Introduction}
\label{Section Introduction}
One of the most classical problems in statistics is the estimation of the parameter vector of a parametrized family of probability distributions. It presents itself in a significant share of applications because parametric models often contribute a reasonable compromise between flexibility in the shape of the statistical model and meaningfulness of the conclusions that can be drawn from the model. As a consequence, all kinds of professions are confronted with the issue of parameter estimation, be it meteorologists, engineers or biologists. Throughout the last decades, a vast amount of highly focused estimation procedures for all kinds of situations have been provided, but the procedure that is arguably used most often remains the maximum likelihood estimator. Apart from its (asymptotic) optimality properties, its popularity is presumably in direct relation with its universality: For the professions mentioned above, and many more, whose prime interest is not the study of sophisticated statistical procedures, it is essential to have at hand a method that is both, easily communicated and applicable to a wide range of model assumptions. A second class of methods incorporates the idea of using as an estimator the value that minimizes some goodness-of-fit measure. To implement this type of estimators, the empirical distribution, quantile or characteristic function is compared to its theoretical counterpart from the underlying parametric model in a suitable distance, and the term is minimized over the parameter space, see \cite{W:1957}, or \cite{P:1981} for an early bibliography. These procedures provide some freedom in adapting the estimation method to the intended inferences from the model and they regularly possess good robustness properties [see \cite{PS:1980} as well as \cite{M:1981}]. An example which was discussed recently, and which goes by the name of minimum CRPM estimation, see \cite{GRWG:2005}, is tailored to the practice of issuing forecasts: As argued by \cite{BGR:2007}, a good probabilistic forecast minimizes a (strictly) proper scoring rule such as the 'CRPM' [\cite{GR:2007}], and after constructing a suitable model it appears somewhat more natural to use as an estimator the one that minimizes the scoring rule instead of a classical estimation method like maximum likelihood [for a comparison see \cite{GMMZ:2018}]. As it happens, these rather universal procedures listed above easily run into computational hardships. Just consider that even for 'basic' models, density functions can take complicated forms, and distribution or characteristic functions, or even normalization constants, may be nowhere near to an explicit formula. This is where we want to tie on. In a recent work, \cite{BE:2018} established distributional characterizations that, from a practical point of view, are comparable to the characterization of a probability distribution through its distribution function. Their results, which are given in terms of the derivative of a density function and the density itself, provide explicit formulae that simplify the dependence of the terms on the parameters (even for rather complicated models), and extend characterizations via the zero-bias- or equilibrium transformation [\cite{GR:1997}, \cite{PR:2011}, respectively] that arise in the context of Stein's method, cf. \cite{CGS:2011}. The aim of this work is to investigate these characterizations, which where already used to construct goodness-of-fit tests [see \cite{BE:2019:1}, \cite{BE:2019:2}], more closely in the context of parameter estimation. An advantage of the resulting estimators lies in the way the density function of the underlying model appears in the characterization, and thus also in the estimation method. When considering for some (positive) density function $p$ the quotient $\tfrac{p^\prime}{p}$, the term no longer depends on the integration constant which ensures that the function integrates to one, but only on the functional form of the density. As indicated before, our estimators depend on the underlying model precisely via this quotient, so they are applicable in cases where the normalization constant is unknown. Models of this type occur (though often in discrete settings) in such applied areas as image modeling [using Markov random fields, see \cite{L:2009}] and machine learning, or in any other area where models are complex enough to render the calculation of the normalization constant impractical. For more specific discussions of such applications, we refer to the introduction of the work by \cite{UKTM:2019}. The problem was already addressed by \cite{H:2005}, who set out to find an estimation method which only takes into account the functional form of a density. The approach introduced there goes by the name of 'score matching', and the estimation method involves terms of the form $\tfrac{p^{\prime\prime}}{p} - \tfrac{1}{2} \big( \tfrac{p^\prime}{p} \big)^2$ and hence does not depend on the normalization constant either. In the univariate case we discuss here, our method provides a good supplement as it contains no second derivatives and may thus be applicable to cases where other methods fail. Also note that several other approaches by \cite{PGH:2010}, \cite{MH:2019}, and \cite{UMK:2019}, are available. Later on we also discuss noise-contrastive estimation, a concept introduced by \cite{GH:2010}. All these references indicate that statistical inference for non-normalized models is a topic of very recent investigation that also interests researcher in machine learning, a fact which we further allude to at the end of the following section. \\

In Section \ref{SEC the new estimators} we introduce this new class of parameter estimators that are comparable, in their range of applicability in the given setting, to the maximum likelihood and minimum Cram\'{e}r-von Mises distance estimators [as discussed by \cite{PS:1980} or \cite{PW:1981}]. We rigorously deal with the existence and measurability of our estimators in Section \ref{SEC Existence and measurability}. In Section \ref{SEC Consistency} we provide results on consistency. Thereafter, we give as (normalized) examples the exponential- (Section \ref{SEC exponential distribution}), the Rayleigh- (Section \ref{SEC Rayleigh distribution}), and the Burr Type XII distribution (Section \ref{SEC Burr distribution}). For each of the three parametric models we compare our new method to classical methods like the maximum likelihood and minimum Cram\'{e}r-von Mises distance estimator in competitive Monte Carlo simulation studies. The Burr distribution [cf. \cite{B:1942}, \cite{R:1977}, \cite{T:1980}, Section 6.2 of \cite{KK:2003}, or \cite{K:2017}] as a model is relevant in econometrics, initiated by \cite{SM:1976} [see also \cite{S:1983}], and other areas like engineering, hydrology, and quality assurance, see \cite{SG:1993} for corresponding references. However, the parameter estimation is non-trivial and can even cause computational issues. Thus, providing a new estimation method could prove useful in applications. In Section \ref{SEC exponential-polynomial models} we discuss an exponential-polynomial model for which the normalization constant is intractable, and we compare the new estimators with the score matching and noise-contrastive estimation approaches.

\section{The new estimators}
\label{SEC the new estimators}

To be specific, recall that the problem of parameter estimation for continuous, univariate probability distributions presents itself as follows. Consider for $\Theta \subset \R^d$ a parametric family of probability density functions
\begin{align*}
\mathfrak{P}_\Theta = \big\{ p_\vartheta \, | \, \vartheta \in \Theta \big\},
\end{align*}
and let $X_1, \dots, X_n$ be a sample consisting of independent real-valued random variables with a distribution from $\mathfrak{P}_\Theta$, that is, there exists some $\vartheta_0 \in \Theta$ such that $X_i$ has density function $p_{\vartheta_0}$ ($X_i \sim p_{\vartheta_0}$, for short) for $i = 1, \dots, n$. Denote with $P_\vartheta$ the distribution function corresponding to $p_\vartheta$. The task is to construct an estimator of the unknown $\vartheta_0$ based on $X_1, \dots, X_n$. 

For the construction of our new estimation method, we first recall in a non-technical fashion a famous distributional characterization that can be traced back to Charles Stein, see Chapter VI of \cite{S:1986}. In the more elaborated version of \cite{LS:2013} it establishes that, given a suitable probability density function $p$, the distribution of a real-valued random variable $X$ is given through the density function at hand if, and only if,
\begin{align*}
\E \left[ f^{\prime}(X) + \frac{p^{\prime}(X)}{p(X)} \, f(X) \right] = 0
\end{align*}
for a large enough class of suitably chosen test functions $f$. Motivated by the well-known zero-bias distribution, \cite{BE:2018} used the above characterization in a recent publication to derive explicit identities which retain the essence of the characterizing property. Indeed, they were able to derive from the Stein characterization that, for a suitable density function $p$ on the positive axis with few technical assumptions (which we adopt below), the distribution of a positive random variable $X$ (satisfying a weak integrability property) is given through $p$ if, and only if, the distribution function $F_X$ corresponding to $X$ satisfies
\begin{align} \label{characterization Betsch Ebner}
F_X(t) = \E \left[ - \frac{p^{\prime}(X)}{p(X)} \, \min\{ X, t \} \right], \quad t > 0.
\end{align}

As we intent to use this result as a foundation for our estimation method in parametric models for non-negative quantities, assume that the support of each density function in $\mathfrak{P}_\Theta$ is $(0, \infty)$. In particular, suppose that each $p_\vartheta$ is positive and continuously differentiable on $(0, \infty)$. Also assume that 
\begin{align*}
\int_{0}^\infty |x| \, \big| p^\prime_\vartheta(x) \big| \, \mathrm{d}x < \infty \quad \text{and} \quad
\sup_{x \, > \, 0} \frac{p^\prime_\vartheta(x) \, \min\{ P_\vartheta(x), 1 - P_\vartheta(x) \}}{{p_\vartheta}^2(x)} < \infty.
\end{align*}
Moreover, suppose that $\lim_{x \, \searrow \, 0} \tfrac{P_\vartheta(x)}{p_\vartheta(x)} = 0$. These presumptions where made by \cite{BE:2018} to derive the characterization given above, and they are straight forward to check for most common density functions. Particularly the last condition is exhaustively discussed in Proposition 3.7 of \cite{D:2015}. Let $X$ be a positive random variable with
\begin{align} \label{Integrability condition for X in the characterization}
\E \left| \frac{p_\vartheta^{\prime}(X)}{p_\vartheta(X)} \, X \right| < \infty, \quad \vartheta \in \Theta,
\end{align}
and define the function
\begin{align*}
\eta(t, \vartheta) = \E\left[ - \frac{p^\prime_\vartheta(X)}{p_\vartheta(X)} \, \min\{ X, t \} \right] - F_{X}(t)
\end{align*}
for $(t, \vartheta) \in (0, \infty) \times \Theta$. Then, the characterization of \cite{BE:2018}, as built up in Equation (\ref{characterization Betsch Ebner}) and as given in their Corollary 3, states that $X$ has density function $p_{\vartheta}$ if, and only if, $\eta(t, \vartheta) = 0$ for every $t > 0$. Therefore, if we assume initially that $X \sim p_{\vartheta_0}$ [note that (\ref{Integrability condition for X in the characterization}) is satisfied by requirement on $p_\vartheta$], then 
\begin{align*}
\big\lVert \eta(\cdot \, , \vartheta) \big\rVert_{L^q} = 0 \quad \text{if, and only if,} \quad \vartheta = \vartheta_0.
\end{align*}
Here, $L^q = L^q\big( (0, \infty), \mathcal{B}(0, \infty), w(t) \, \mathrm{d}t \big)$, $1 \leq q < \infty$, denote the usual $L^q$-spaces over $(0, \infty)$, $w$ is a positive and integrable weight function, and for $f \in L^q$, $g \in L^{q^\prime}$ ($1/q + 1/q^{\prime} = 1$)
\begin{align*}
\lVert f \rVert_{L^q} = \left(\int_0^\infty |f(t)|^q \, w(t) \, \mathrm{d}t\right)^{1/q}, \quad \langle f, g \rangle_{L^q} = \int_0^\infty f(t) \, g(t) \, w(t) \, \mathrm{d}t
\end{align*}
are the usual norm and duality in $L^q$. Thus, with an empirical version
\begin{align} \label{Empirical version of eta}
\eta_n(t, \vartheta) 
= - \frac{1}{n} \sum_{j=1}^{n} \frac{p^\prime_\vartheta(X_j)}{p_\vartheta(X_j)} \, \min\{ X_j, t \} - \frac{1}{n} \sum_{j=1}^{n} \mathds{1}\{ X_j \leq t \}
\end{align}
of $\eta$, based on a sample of independent and identically distributed (i.i.d.) random variables $X_1, \dots, X_n$ with $X_1 \sim p_{\vartheta_0}$, a reasonable estimator for the unknown $\vartheta_0$ is
\begin{align} \label{Estimator of theta}
\widehat{\vartheta}_{n, q}
= \arg\min\big\{ \lVert \eta_n(\cdot \, , \vartheta) \rVert_{L^q} \, | \, \vartheta \in \Theta \big\} \,
\Big( = \arg\min\big\{ \lVert \eta_n(\cdot \, , \vartheta) \rVert_{L^q}^q \, | \, \vartheta \in \Theta \big\} \Big),
\end{align}
that is, we choose $\widehat{\vartheta}_{n, q}$ such that $\lVert \eta_n(\cdot \, , \widehat{\vartheta}_{n, q}) \rVert_{L^q} \leq \lVert \eta_n(\cdot \, , \vartheta) \rVert_{L^q}$ for each $\vartheta \in \Theta$. Heuristically, $\lVert \eta_n(\cdot \, , \vartheta) \rVert_{L^q}$ approximates $\lVert \eta(\cdot \, , \vartheta) \rVert_{L^q}$, so $\widehat{\vartheta}_{n, q}$ should provide an estimate for the minimum of $\vartheta \mapsto \lVert \eta(\cdot \, , \vartheta) \rVert_{L^q}$ which coincides with $\vartheta_0$, the (unique) zero of this function. At this point of course, there arise questions of existence and measurability of such an estimator, and we will handle these questions in full detail in Section \ref{SEC Existence and measurability}. Intuitively, one might argue to replace $F_X$ and the empirical distribution function in the definition of $\eta$ and $\eta_n$, respectively, with the theoretical distribution function $P_{\vartheta}$. However, there is a bit of a technical point involved, and the characterizations by \cite{BE:2018} do not include results that give a rigorous handle for this slightly (yet decisively) different situation. There are, however, similar characterizations for univariate distribution with other supports than the positive half axis. We allude to that setting in Section \ref{SEC notes and comments}. Note that the availability of the term $\tfrac{p^\prime_\vartheta}{p_\vartheta}$ for the model in consideration is rather essential. If this term is not amenable explicitly, it might still be calculable using numerical differentiation (and so $\eta_n$ and the estimator could be calculated numerically), but it would make it hard to theoretically justify the validity of the conditions on $p_\vartheta$. In our experience, however, the term $\tfrac{p^\prime_\vartheta}{p_\vartheta}$ is readily available whenever $p_\vartheta$ can be differentiated explicitly, and this seems a manageable assumption. \\

As we have outlined above, our new estimators are eventually based on Stein characterizations which rely on some suitable class of test functions [for an overview in the univariate case, and a record of the vast amount of literature on these identities, see \cite{LRS:2017}]. The goal of \cite{BE:2018} was to derive from these Stein identities new characterizations that no longer involve the classes of test functions. While this approach apparently leads to feasible applications in statistics, other methods are based directly on the Stein characterizations, using Stein discrepancies which gradually become popular in machine learning. The idea in the context of parameter estimation, in heuristic terms, boils down to choosing as a parameter estimator the value which (approximately) minimizes
\begin{align*}
\sup_{f} \left| \E \left[ f^{\prime}(X) + \frac{p_\vartheta^{\prime}(X)}{p_\vartheta(X)} \, f(X) \right] \right|,
\end{align*}
where the supremum is over all test functions in consideration. By the Stein characterization detailed above, the expectation is $0$ for every test function precisely when $\vartheta = \vartheta_0$, as we assume that $X \sim p_{\vartheta_0}$. However, it is not clear how to calculate the supremum in practice taking that the class of test functions is very large. The theory developed around Stein discrepancies has produced different formal methods to evaluate such terms. Other than the fact that they are based on the Stein characterization, the identities derived by \cite{BE:2018} are not related to the framework of Stein discrepancies, and so it is surprising that merely measuring the difference between the functions in (\ref{characterization Betsch Ebner}) in an $L^q$-norm, which is what we do to construct our estimators, leads back to so-called feature Stein discrepancies. Indeed, upon defining the 'feature' function $\Phi(x, t) = \min\{ x, t \}$, $x, t > 0$, and considering the Langevin-Stein operators
\begin{align*}
(\mathcal{T}_\vartheta f)(x, t) = p_\vartheta(x)^{-1} \, \partial_x \big( p_\vartheta(x) \cdot f(x, t) \big)
\end{align*}
as applied to suitable functions $f : (0, \infty)^2 \to \R$, we obtain
\begin{align*}
\lVert \eta(\cdot \, , \vartheta) \rVert_{L^q}^2
= \left\lVert \E\left[ - \frac{p^\prime_\vartheta(X)}{p_\vartheta(X)} \, \min\{ X, \cdot \} \right] - F_{X}(\cdot) \right\rVert_{L^q}^2
= \left\lVert \int_0^\infty (\mathcal{T}_\vartheta \Phi)(x, \cdot) \, \mathrm{d}P_{\vartheta_0}(x) \right\rVert_{L^q}^2,
\end{align*}
which is the right-hand side of Equation (1) in the paper by \cite{HM:2018}. So by retracing their calculation,
\begin{align*}
\lVert \eta(\cdot \, , \vartheta) \rVert_{L^q}^2 
= \Phi\mathrm{SD}^2_{\Phi, q}(P_{\vartheta_0}, P_\vartheta) 
= \sup_{g \in \mathcal{G}_{\Phi, q}} \left| \E \left[ \frac{\partial_x \big( g(X) \cdot p_\vartheta(X) \big)}{p_\vartheta(X)} \right] \right|^2 ,
\end{align*}
where $\mathcal{G}_{\Phi, q}$ is the class of test functions as defined by \cite{HM:2018} (the precise form of which is inessential at this point). This means that we can embed our setting into the framework of these feature Stein discrepancies, as the construction of our estimator cumulating in (\ref{Estimator of theta}) corresponds to minimizing the quantity at the beginning of this paragraph which sought to motivate these discrepancies. Now, of course, the starting point of our estimation method being the characterizations by \cite{BE:2018}, we already had our estimator at hand explicitly and could choose the feature function accordingly. Still, the fact that both the characterization of \cite{BE:2018} and the (feature) Stein discrepancy approach (for the above feature function), when translated into an estimation method, lead to the same procedure is remarkable and deems it worthwhile to study the method further, as we were assured that it can be rather hard to find explicit examples for which the Stein discrepancy approach is feasible in practice. 

To complete this insightful tour into the realm of Stein discrepancies, we mention some contributions of various solutions to statistical problems based on those discrepancies. In particular, \cite{CSG:2016}, \cite{LLJ:2016}, and \cite{YLRN:2018} construct tests of fit, \cite{GM:2015} measure sample quality, and \cite{BBDGM:2019} develop estimation methods for non-normalized statistical models.

\section{Existence and measurability}
\label{SEC Existence and measurability}
We discuss the measurability properties of $\eta_n$ and derive an existence result for a measurable version of (approximate) estimators of the type in (\ref{Estimator of theta}). The result that is central to us in this section can be found in Chapter III of \cite{CV:1977} [see the references therein and Chapter 8 by \cite{C:2013} for further background]. Before we summarize these results, recall that a Suslin space is a Hausdorff topological space which is the image of a separable, completely metrizable topological space under a continuous map [for an overview, consult Chapter II of \cite{S:1973}]. See also Remark \ref{RMK comments to selection theorem} in Appendix \ref{Appendix existence section} for more information.
\begin{theorem} \label{THM measurable selection}
	Let $(\Omega, \mathcal{F}, \mathbb{P})$ be a complete probability space and $(\mathfrak{S}, \mathcal{O}_{\mathfrak{S}})$ a Suslin topological space with Borel-$\sigma$-field $\mathcal{B}(\mathfrak{S})$. Assume that $\Gamma$ maps $\Omega$ into the non-empty subsets of $\mathfrak{S}$, and that
	\begin{align*}
	\mathrm{graph}(\Gamma)
	= \big\{ (\omega, x) \in \Omega \times \mathfrak{S} \, | \, x \in \Gamma(\omega) \big\}
	\in \mathcal{F} \otimes \mathcal{B}(\mathfrak{S}).
	\end{align*}
	Then, there exists an $\big(\mathcal{F}, \mathcal{B}(\mathfrak{S})\big)$-measurable map $\widehat{\vartheta} : \Omega \to \mathfrak{S}$ such that $\widehat{\vartheta}(\omega) \in \Gamma(\omega)$ for every $\omega \in \Omega$. Additionally, if $\psi : \Omega \times \mathfrak{S} \to \overline{\R}$ is $\big(\mathcal{F} \otimes \mathcal{B}(\mathfrak{S}), \overline{\mathcal{B}}\big)$-measurable, then
	\begin{align*}
	m(\omega) = \inf_{x \, \in \, \Gamma(\omega)} \psi(\omega, x) \quad \text{and} \quad
	M(\omega) = \sup_{x \, \in \, \Gamma(\omega)} \psi(\omega, x)
	\end{align*}
	are $(\mathcal{F}, \overline{\mathcal{B}})$-measurable. Here, $(\overline{\R}, \overline{\mathcal{B}})$ denotes the extended real line with its usual $\sigma$-field, and we write $\otimes$ for the product of $\sigma$-fields.
\end{theorem}

To apply Theorem \ref{THM measurable selection}, we first have to investigate the measurability properties of $\eta_n$. In the setting of Section \ref{SEC the new estimators}, assume the following regularity condition.
\begin{itemize}
	\item[\bf(R1)]{The map $\Theta \ni \vartheta \mapsto \frac{p^{\prime}_\vartheta(x)}{p_\vartheta(x)}$ is continuous for every $x > 0$.}
\end{itemize}
Let $(\Omega, \mathcal{F}, \mathbb{P})$ be a complete probability space, which is assumed to underlie all random quantities of the previous and subsequent sections. Notice that the function $\eta_n$ defined in (\ref{Empirical version of eta}) depends on the random variables $X_1, \dots, X_n$ defined on $\Omega$, hence $\eta_n$ (as a random quantity) can be understood as a map $\Omega \times (0, \infty) \times \Theta \to \R$. Exploiting the structure of $\eta_n$, we obtain the following lemma. The proof is simple, and the basic thoughts can be found in Appendix \ref{Appendix existence section}.
\begin{lemma} \label{LEMMA measurability eta_n}
	The map $\eta_n : \Omega \times (0, \infty) \times \Theta \to \R$ from (\ref{Empirical version of eta}) is $\big( \mathcal{F} \otimes \mathcal{B}(0, \infty) \otimes \mathcal{B}(\Theta), \mathcal{B}^1 \big)$-measurable. Moreover, as an element in $L^q$, $\eta_n : \Omega \times \Theta \to L^q$ is $\big( \mathcal{F} \otimes \mathcal{B}(\Theta), \mathcal{B}(L^q) \big)$-measurable. In particular,
	\begin{align*}
	(\omega, \vartheta) \mapsto \big\lVert \eta_n(\omega, \cdot \, , \vartheta) \big\rVert_{L^q}
	\end{align*}
	is an $\big( \mathcal{F} \otimes \mathcal{B}(\Theta), \mathcal{B}[0, \infty) \big)$-measurable mapping.
\end{lemma}
Similar measurability results hold for $\eta : (0, \infty) \times \Theta \to \R$. For the remainder of this work assume that
\begin{itemize}
	\item[\bf(R2)]{the parameter space $\emptyset \neq \Theta \subset \R^d$ is a Borel set in $\R^d$}.
\end{itemize}
As such, $\Theta$ is a Suslin topological space [see Proposition 8.2.10 from \cite{C:2013}] with the subspace topology induced by $\R^d$. It is also a metric space with the standard metric in $\R^d$ restricted to $\Theta$. For $n \in \N$, let $\varepsilon_n$ be positive random variables such that $\varepsilon_n \to 0$ $\mathbb{P}$-almost surely (a.s.), as $n \to \infty$. Define $\psi_{n, q}(\omega, \vartheta) = \big\lVert \eta_n(\omega, \cdot \, , \vartheta) \big\rVert_{L^q}$ which, by Lemma \ref{LEMMA measurability eta_n}, is a product-measurable function from $\Omega \times \Theta$ into $[0, \infty)$. Theorem \ref{THM measurable selection} implies that $m_{n, q}(\omega) = \inf_{\vartheta \, \in \, \Theta} \psi_{n, q}(\omega, \vartheta)$ is $(\mathcal{F}, \mathcal{B}^1)$-measurable. Hence the set-valued function
\begin{align} \label{set of approximate estimators}
\Gamma_{n, q}(\omega)
= \Big\{ \vartheta \in \Theta \, \Big| \, \psi_{n ,q}(\omega, \vartheta) \leq m_{n, q}(\omega) + \varepsilon_n(\omega) \Big\}
\end{align}
has a measurable graph. By construction, $\Gamma_{n, q}$ takes as values only non-empty subsets of $\Theta$. In fact, $\Gamma_{n, q}(\omega)$  is also closed in $\Theta$ for every $\omega \in \Omega$, see Remark \ref{RMK closedness of Gamma_n,q} in Appendix \ref{Appendix existence section}. Theorem \ref{THM measurable selection} yields the existence of an $\big(\mathcal{F}, \mathcal{B}(\Theta)\big)$-measurable map $\widehat{\vartheta}_{n, q} : \Omega \to \Theta$ with $\widehat{\vartheta}_{n, q}(\omega) \in \Gamma_{n, q}(\omega)$, which is, by definition of $\Gamma_{n, q}$,
\begin{align} \label{Estimator, formal defintion}
\Big\lVert \eta_n\big(\omega, \cdot \, , \widehat{\vartheta}_{n, q}(\omega)\big) \Big\rVert_{L^q}
\leq \inf\limits_{\vartheta \, \in \, \Theta} \Big\lVert \eta_n(\omega, \cdot \, , \vartheta) \Big\rVert_{L^q} + \varepsilon_n(\omega)
\end{align}
for each $\omega \in \Omega$ or, equivalently, 
\begin{align*}
\psi_{n, q} \Big( \omega, \widehat{\vartheta}_{n, q}(\omega) \Big)
\leq \inf\limits_{\vartheta \, \in \, \Theta} \psi_{n, q}\big( \omega, \vartheta \big) + \varepsilon_n(\omega) .
\end{align*}
Whenever we refer to an estimator that satisfies (\ref{Estimator of theta}), we mean precisely such an (approximate) measurable version. This settles the existence problem and for our asymptotic studies we have measurability of $\widehat{\vartheta}_{n, q}$ at hand.

\section{Consistency}
\label{SEC Consistency}
In this section, we investigate the asymptotic behavior of our estimators. Unfortunately, we cannot apply the general results for minimum distance estimators given by \cite{M:1984}, since a major assumption in that work is that the term in the norm is differentiable (with respect to $\vartheta$) with derivative not depending on $\omega$, that is, in a sense, the parameter and the 'uncertainty' have to be separated, which is clearly not the case in our setting. Thus, we need to deal with the empirical process involved.

Assume the setting from Section \ref{SEC the new estimators}. For brevity, we keep the notation $\psi_{n, q}(\vartheta) \big(= \psi_{n, q}(\omega, \vartheta)\big) = \big\lVert \eta_n(\omega, \cdot \, , \vartheta) \big\rVert_{L^q}$ and set $\psi_q(\vartheta) = \big\lVert \eta(\cdot \, , \vartheta) \big\rVert_{L^q}$. Recall from the construction that $\widehat{\vartheta}_{n, q}$ (approximately) minimizes $\psi_{n, q}$ [see (\ref{Estimator, formal defintion})], and $\vartheta_0$ is the unique minimum of $\psi_q$. The heuristic of the consistency statement proven in this section is as follows. If $\psi_{n, q}$ converges to $\psi_q$ in a suitable function space, then the random minimal points $\widehat{\vartheta}_{n, q}$ converge to $\vartheta_0$. In order to establish convergence of $\psi_{n, q}$, we need the functions to be sufficiently smooth in $\vartheta$. In most applications the mapping $\vartheta \mapsto \frac{p_\vartheta^\prime(x)}{p_\vartheta(x)}$ will be continuously differentiable for every $x > 0$, which can often be used to derive the following regularity condition.
\begin{itemize}
	\item[\bf(R3)]{For each non-empty compact subset $K$ of $\Theta$ there exists some $0 < \alpha = \alpha_K < \infty$ and a measurable function $H = H_K : (0, \infty) \to [0, \infty)$ with $\E \big[ H(X) \, X \big] < \infty$ such that
		\begin{align*}
		\left| \frac{p_{\vartheta^{(2)}}^\prime(x)}{p_{\vartheta^{(2)}}(x)} - \frac{p_{\vartheta^{(1)}}^\prime(x)}{p_{\vartheta^{(1)}}(x)} \right|
		\leq H(x) \, \big| \vartheta^{(2)} - \vartheta^{(1)} \big|^\alpha,
		\end{align*}
		for every $x > 0$ and all $\vartheta^{(1)}, \vartheta^{(2)} \in K$.}
\end{itemize}
Now, let $K \neq \emptyset$ be an arbitrary compact subset of $\Theta$. Then on $\Omega$ and for $\vartheta^{(1)}, \vartheta^{(2)} \in K$, we have
\begin{align*}
\Big| \psi_{n, q}\big( \vartheta^{(2)} \big) - \psi_{n, q}\big( \vartheta^{(1)} \big) \Big|
\leq \big| \vartheta^{(2)} - \vartheta^{(1)} \big|^\alpha \cdot \left(\int_0^\infty w(t) \, \mathrm{d}t\right)^{1/q} \cdot \frac{1}{n} \sum_{j=1}^{n} H(X_j) \, X_j
\end{align*}
with $H$ and $\alpha$ as in (R3). In particular, $K \ni \vartheta \mapsto \psi_{n, q}(\omega, \vartheta)$ is continuous for every $\omega \in \Omega$, and, by Lemma \ref{LEMMA measurability eta_n}, it constitutes a product measurable map. This already implies that $\vartheta \mapsto \psi_{n, q}(\vartheta)$ is a random element of $C(K)^+$ [see Lemma 3.1 of \cite{K:2002}], the space of continuous functions from $K$ to $[0, \infty)$ which is a complete, separable metric space (endowed with the usual metric that induces the uniform topology). From (R3) it also follows that $K \ni \vartheta \mapsto \psi_q(\vartheta)$ is an element of $C(K)^+$. We can now state the convergence results for $\psi_{n, q}$ that are essential for our consistency proof.
\begin{lemma} \label{LEMMA convergence of psi_n,q}
	In the setting of Section \ref{SEC the new estimators}, assume that (R1) - (R3) are satisfied. Let $K \neq \emptyset$ be a compact subset of $\Theta$. Then $\psi_{n, q} \longrightarrow \psi_q$ in $C(K)^+$ $\mathbb{P}$-a.s., as $n \to \infty$. Moreover,
	\begin{align*}
	\inf_{\vartheta \, \in \, F} \big\lVert \eta_n(\, \cdot \, , \vartheta) \big\rVert_{L^q}
	= \inf_{\vartheta \, \in \, F} \psi_{n, q}(\vartheta)
	\longrightarrow \inf_{\vartheta \, \in \, F} \psi_q(\vartheta)
	= \inf_{\vartheta \, \in \, F} \big\lVert \eta(\, \cdot \, , \vartheta) \big\rVert_{L^q}
	\end{align*}
	$\mathbb{P}$-a.s., as $n \to \infty$, for every non-empty closed subset $F$ of $K$.
\end{lemma}
The proof of this lemma is rather technical and deferred to Appendix \ref{Appendix consistency section}. Note that the term $\inf_{\vartheta \, \in \, F} \psi_{n, q}(\vartheta)$ is a random variable by Theorem \ref{THM measurable selection} (cf. the measurability of $m_{n, q}$ in the previous section). The following theorem uses the above lemma to establish consistency. In the second statement, we assume that the parameter space $\Theta$ is compact, thus rendering Lemma \ref{LEMMA convergence of psi_n,q} applicable on the whole of $\Theta$, which will turn out essential to prove strong consistency. For most practical purposes this is sufficient, when parameters relevant for modeling in applications can be taken to stem from some (huge) compact set. Note that with this compactness assumption we actually do not need the $\varepsilon_n$-term in (\ref{Estimator, formal defintion}) since $\psi_{n, q}$ is lower semi-continuous by (R1) and Fatou's lemma, and thus attains its minimum in $\Theta$. The first statement of the following theorem shows that if the sequence $\widehat{\vartheta}_{n, q}$ is already known to be tight, no compactness assumption is needed, but we can only expect weak consistency in general, thus denoting by '$\stackrel{\mathbb{P}}{\longrightarrow}$' convergence in probability. After the proof of the theorem, we provide an insight in which cases this is possible (Remark \ref{RMK consistency in convex case}).
\begin{theorem}{[Consistency]} \label{THM consistency}
	Take the setting from Section \ref{SEC the new estimators}, let $\psi_{n, q}$, $\psi_q$ be as above, and consider $\widehat{\vartheta}_{n, q}$ from (\ref{Estimator, formal defintion}). Further assume that (R1) -- (R3) are satisfied. 
	\begin{itemize}
		\item[(i)] {If $\big\{ \widehat{\vartheta}_{n, q} \big\}_{n \, \in \, \N}$ is tight in $\Theta$, then $\widehat{\vartheta}_{n, q} \stackrel{\mathbb{P}}{\longrightarrow} \vartheta_0$, as $n \to \infty$.}
		\item[(ii)] {If $\Theta$ is compact, then $\widehat{\vartheta}_{n, q} \longrightarrow \vartheta_0$ $\mathbb{P}$-a.s., as $n \to \infty$.}
	\end{itemize}
\end{theorem}
\begin{prf}
	In the proof of $(i)$ we follow Theorem 3.2.2 from \cite{VW:2000}, but we adapt the reasoning to our setting, using the measurability properties we established in Section \ref{SEC Existence and measurability}, and Lemma \ref{LEMMA convergence of psi_n,q}. For completeness, as well as to prepare the proof of the second result, we give a full proof. We start with a preliminary observation, establishing that the minimum at $\vartheta_0$ is (locally) well separated. If $K$ is a compact subset of $\Theta$, and $O$ an open subset of $\R^d$ which contains $\vartheta_0$, then
	\begin{align} \label{minimum is well separated}
	0 = \psi_q(\vartheta_0) < \inf_{\vartheta \, \in \, K \setminus O} \psi_q(\vartheta) .
	\end{align} 
	Indeed, if this is not the case, we find a sequence $\vartheta^{(k)} \in K \setminus O$ such that $\psi_q\big(\vartheta^{(k)}\big) \longrightarrow 0$ as $k \to \infty$. Since $K \setminus O$ is compact, there exists a subsequence $\{ \vartheta^{(k_j)} \}_{j \, \in \, \N}$ and some $\vartheta^* \in K \setminus O$ such that $\vartheta^{(k_j)} \longrightarrow \vartheta^*$ as $j \to \infty$. By continuity of $\psi_q$, it holds that $\psi_q(\vartheta^*) = \lim_{j \, \to \, \infty} \psi_q\big(\vartheta^{(k_j)}\big) = 0$, but $K \setminus O \ni \vartheta^* \neq \vartheta_0 \in O$ which is a contradiction to the fact that $\vartheta_0$ is the unique zero of $\psi_q$. \\
	
	Now, let $\varepsilon, \delta > 0$. Choose a compact subset $K = K_\delta \subset \Theta$ with $\sup_{n \, \in \, \N} \mathbb{P}\big( \widehat{\vartheta}_{n, q} \notin K \big) < \delta$, and define $F = K \setminus B_{\varepsilon}(\vartheta_0)$, where $B_{\varepsilon}(\vartheta_0)$ denotes the open ball in $\R^d$ of radius $\varepsilon$ around $\vartheta_0$. Applying Lemma \ref{LEMMA convergence of psi_n,q} and (\ref{minimum is well separated}) to $K$ and $F$, together with (\ref{Estimator, formal defintion}) and the Portmanteau theorem [cf. Theorem 2.1 of \cite{B:1968}], we get
	\begin{align*}
	\limsup_{n \, \to \, \infty} \mathbb{P} \Big( \big| \widehat{\vartheta}_{n, q} - \vartheta_0 \big| \geq \varepsilon \Big)
	&\leq \limsup_{n \, \to \, \infty} \mathbb{P} \Big( \widehat{\vartheta}_{n, q} \in F \Big) + \limsup_{n \, \to \, \infty} \mathbb{P} \Big( \widehat{\vartheta}_{n, q} \notin K \Big) \\
	&\leq \limsup_{n \, \to \, \infty} \mathbb{P} \Big( \inf_{\vartheta \, \in \, F} \psi_{n, q}(\vartheta) \leq \psi_{n, q}(\vartheta_0) + \varepsilon_n \Big) + \delta \\
	&\leq \mathbb{P} \Big( \inf_{\vartheta \, \in \, F} \psi_q(\vartheta) \leq \psi_q(\vartheta_0) \Big) + \delta \\
	&= \delta.
	\end{align*}
	Note that if $F = \emptyset$, the inequality holds trivially. Since both $\varepsilon$ and $\delta$ were arbitrary, the claim follows. For this first part of the proof, we only needed the convergences provided by Lemma \ref{LEMMA convergence of psi_n,q} to be valid in probability. For the following proof of $(ii)$, we rely on the stronger result. The arguments we use are scattered over Section 3 of the work by \cite{S:1970}. For reasons alluded to in Remark \ref{RMK comments to selection theorem}, and since that work contains some typos, we provide the adapted arguments. Let $\varepsilon > 0$ and define $\beta_\varepsilon = \inf_{\vartheta \, \in \, \Theta \setminus B_\varepsilon(\vartheta_0)} \psi_q(\vartheta)$. By (\ref{minimum is well separated}), we have $\beta_\varepsilon > 0$. Using the well-known equivalent criterion for almost sure convergence, Lemma \ref{LEMMA convergence of psi_n,q} gives
	\begin{align*}
	\lim_{n \, \to \, \infty} \mathbb{P}\left( \bigcup_{k \geq n} \left\{ \sup_{\vartheta \, \in \, \Theta \setminus B_\varepsilon(\vartheta_0)} \Big| \psi_{k, q}(\vartheta) - \psi_q(\vartheta) \Big| \geq \frac{\beta_\varepsilon}{2} \right\} \right)
	= 0.
	\end{align*}
	By definition of $\beta_\varepsilon$ this implies
	\begin{align*}
	\lim_{n \, \to \, \infty} \mathbb{P}\left( \bigcup_{k \geq n} \left\{ \inf_{\vartheta \, \in \, \Theta \setminus B_\varepsilon(\vartheta_0)} \psi_{k, q}(\vartheta) \leq \frac{\beta_\varepsilon}{2} \right\} \right)
	= 0.
	\end{align*}
	Moreover, $\psi_{n, q}(\vartheta_0) + \varepsilon_n \longrightarrow \psi_q(\vartheta_0) = 0$ $\mathbb{P}$-a.s., as $n \to \infty$, and thus
	\begin{align*}
	\lim_{n \, \to \, \infty} \mathbb{P}\left( \bigcup_{k \geq n} \left\{ \Big| \psi_{k, q}(\vartheta_0) + \varepsilon_k \Big| \geq \frac{\beta_\varepsilon}{2} \right\} \right)
	= 0.
	\end{align*}
	Putting everything together,
	\begin{align*}
	\limsup_{n \, \to \, \infty} \mathbb{P} \left( \bigcup_{k \geq n} \Big\{ \big| \widehat{\vartheta}_{k, q} - \vartheta_0 \big| \geq \varepsilon \Big\} \right) 
	&\leq \limsup_{n \, \to \, \infty} \mathbb{P} \left( \bigcup_{k \geq n} \left\{ \inf_{\vartheta \, \in \, \Theta \setminus B_\varepsilon(\vartheta_0)} \psi_{k, q}(\vartheta) \leq \psi_{k, q}(\vartheta_0) + \varepsilon_k \right\} \right) \\
	&\leq \limsup_{n \, \to \, \infty} \mathbb{P} \left( \bigcup_{k \geq n} \left\{ \inf_{\vartheta \, \in \, \Theta \setminus B_\varepsilon(\vartheta_0)} \psi_{k, q}(\vartheta) \leq \frac{\beta_\varepsilon}{2} \right\} \right) \\
	&\quad + \limsup_{n \, \to \, \infty} \mathbb{P} \left( \bigcup_{k \geq n} \left\{ \psi_{k, q}(\vartheta_0) + \varepsilon_k \geq \frac{\beta_\varepsilon}{2} \right\} \right) \\
	&= 0,
	\end{align*}
	that is, $\widehat{\vartheta}_{n, q} \longrightarrow \vartheta_0$ $\mathbb{P}$-a.s., as $n \to \infty$.
\end{prf}
\begin{remark}{[A priori tightness of the sequence of estimators]} \label{RMK consistency in convex case}
	We provide a tool for proving tightness of the estimators before having established consistency, which we can use in Theorem \ref{THM consistency} to get consistency even for unbounded parameter spaces. The statement essentially yields that if $\psi_{n, q}$ is strictly convex, $\big\{ \widehat{\vartheta}_{n, q} \big\}_{n \, \in \, \N}$ is tight. More precisely, suppose that conditions (R1) -- (R3) hold. Let $\Theta$ be convex with $\vartheta_0 \in \Theta^{\circ}$, the interior of $\Theta$. Further, let $\psi_{n, q}$ be strictly convex (almost surely). Then the sequence of estimators $\widehat{\vartheta}_{n, q}$ is tight in $\Theta$. The proof is straight-forward and some hints are given in exercise problem 4 in Section 3.2 of \cite{VW:2000} (for more details, find the proof in Appendix \ref{Appendix consistency section}).
\end{remark}

\section{Example: The exponential distribution}
\label{SEC exponential distribution}

Let $\Theta = (0, \infty)$ and $p_\vartheta(x) = \vartheta \exp(- \vartheta x)$, $x > 0$. This trivially is an admissible class of density functions. Moreover, let $\vartheta_0 \in \Theta$, $X \sim p_{\vartheta_0}$, and take a sample $X_1, \dots, X_n$ of i.i.d. copies of $X$. An easy calculation gives
\begin{align*}
\big(\psi_q(\vartheta)\big)^q
&= \int_0^\infty \Big| \E\big[ \vartheta \min\{ X, t \} \big] - \big( 1 - \exp(- \vartheta_0 t) \big) \Big|^q \, w(t) \, \mathrm{d}t \\
&= \Big| \frac{\vartheta}{\vartheta_0} - 1 \Big|^q \int_0^\infty \Big| 1 - \exp(- \vartheta_0 t) \Big|^q \, w(t) \, \mathrm{d}t,
\end{align*}
which nicely illustrates that $\vartheta_0$ is indeed the unique zero of this functions. For the particular choice of weight $w(t) = \exp(- a t)$, $t > 0$, with some tuning parameter $a > 0$, and in the case $q = 2$, straight-forward calculations give
\begin{align*}
\big(\psi_{n, 2}(\vartheta)\big)^2
= \vartheta^2 \Psi_n^{(1)} + \vartheta \Psi_n^{(2)} + \Psi_n^{(3)},
\end{align*}
where
\begin{align*}
\Psi_n^{(1)}
&= \frac{2}{a^3} + \frac{2}{a^2 n^2} \sum_{j=1}^{n} e^{-a X_{(j)}} \Big( X_{(j)} (-n + j - 1) - \frac{1}{a} (2n - 2j + 1) \Big) - \frac{2}{a^2 n^2} \sum_{1 \leq j < k \leq n} X_{(j)} \, e^{-a X_{(k)}} , \\
\Psi_n^{(2)}
&= \frac{2}{a n^2} \sum_{j=1}^{n} e^{-a X_{(j)}} \Big( X_{(j)} (-n + j - 1) - \frac{1}{a} (n - 2j + 1) \Big) - \frac{2}{a n^2} \sum_{1 \leq j < k \leq n} X_{(j)} \, e^{-a X_{(k)}} , \\
\Psi_n^{(3)}
&= \frac{1}{a n^2} \sum_{j=1}^{n} e^{-a X_{(j)}} (2j - 1) ,
\end{align*}
and $X_{(1)} < \dotso < X_{(n)}$ is the ordered sample. Using that $e^{-a X_{(k)}} < e^{-a X_{(j)}}$ $\mathbb{P}$-a.s. for $j < k$, we obtain
\begin{align*}
\Psi_n^{(1)}
\geq \frac{2}{a^3 n^2} \sum_{j=1}^{n} (2n - 2j + 1) \Big( 1 - \big( 1 + a X_{(j)} \big) e^{-a X_{(j)}} \Big),
\end{align*}
and since $1 + a X_{(j)} < e^{a X_{(j)}}$ $\mathbb{P}$-a.s., we have $\Psi_n^{(1)} > 0$ almost surely. Therefore, $\big(\psi_{n, 2}\big)^2$ is strictly convex (almost surely), and has a unique minimum. By Remark \ref{RMK consistency in convex case} and Theorem \ref{THM consistency} $(i)$, the estimator
\begin{align*}
\widehat{\vartheta}_{n, 2}^{(a)}
= \mathrm{arg \, min} \big\{ \psi_{n, 2}(\vartheta) \, \big| \, \vartheta > 0 \big\}
&= \mathrm{arg \, min} \big\{ \big(\psi_{n, 2}(\vartheta)\big)^2 \, \big| \, \vartheta > 0 \big\} \\
&= \mathrm{arg \, min} \big\{ \vartheta^2 \Psi_n^{(1)} + \vartheta \Psi_n^{(2)} + \Psi_n^{(3)} \, \big| \, \vartheta > 0 \big\}
\end{align*}
is consistent for $\vartheta_0$ (over the whole of $\Theta$). Note that we have not made the dependence of $\Psi_n^{(1)}$, $\Psi_n^{(2)}$, and $\Psi_n^{(3)}$ on '$a$' explicit to prevent overloading the notation. With a similar argument as above, we may show that $\Psi_n^{(2)} < 0$ almost surely, thus we can calculate $\widehat{\vartheta}_{n, 2}^{(a)}$ explicitly as
\begin{align*}
\widehat{\vartheta}_{n, 2}^{(a)}
= - \frac{\Psi_n^{(2)}}{2 \Psi_n^{(1)}} .
\end{align*}
To provide insight on the performance of this estimator, we compare it with the maximum likelihood estimator and the minimizer of the mean squared error (for $n \geq 3$) which are given as
\begin{align*}
\widehat{\vartheta}_n^{ML} = \left( \frac{1}{n} \sum_{j = 1}^{n} X_j \right)^{-1} \quad \text{and} \quad
\widehat{\vartheta}_n^{MSE} = \left( \frac{1}{n - 2} \sum_{j = 1}^{n} X_j \right)^{-1},
\end{align*}
respectively, as well as with the minimum Cram\'{e}r-von Mises distance estimator discussed in the introduction, namely
\begin{align*}
\widehat{\vartheta}_n^{CvM}
&= \arg\min \left\{ \int_0^\infty \left( \widehat{F}_n(t) - P_\vartheta(t) \right)^2 \, \mathrm{d}P_\vartheta (t) ~ \Big| ~ \vartheta > 0 \right\} \\
&= \arg\min \left\{ \frac{1}{n} \sum_{j = 1}^{n} \left[ \exp\Big( - 2 \vartheta X_{(j)} \Big) + \exp\Big( - \vartheta X_{(j)} \Big) \cdot \Big( \frac{2j - 1}{n} - 2 \Big) \right] ~ \Big| ~ \vartheta > 0 \right\} ,
\end{align*}
where $P_\vartheta(x) = 1 - \exp(- \vartheta x)$, $x > 0$, denotes the distribution function of the exponential distribution, and where $\widehat{F}_n$ is the empirical distribution function of $X_1, \dots, X_n$. For this comparison we simulate (for fixed values of $n$ and $\vartheta_0$) $D = 100,000$ samples of size $n$ from an exponential distribution with parameter $\vartheta_0$, calculate the values of the estimator for each sample yielding values $\widehat{\vartheta}_1, \dots, \widehat{\vartheta}_D$, and approximate the bias and mean squared error (MSE) via
\begin{align*}
\frac{1}{D} \sum_{k = 1}^{D} \Big( \widehat{\vartheta}_k - \vartheta_0 \Big) \quad \text{and} \quad
\frac{1}{D} \sum_{k = 1}^{D} \Big( \widehat{\vartheta}_k - \vartheta_0 \Big)^2
\end{align*}
for each of the above estimators. We perform all simulations with Python 3.7.2 (as provided by the Python Software Foundation, \href{https://www.python.org}{https://www.python.org}, accessed 28 August 2019). For the minimization required to calculate the minimum Cram\'{e}r-von Mises distance estimator, we choose as initial value the maximum likelihood estimator and use a sequential least squares programming method ('SLSQP') [cf. \cite{K:1988}] implemented in the 'optimize.minimize' function of the Python module 'scipy', see \cite{JOP:2001}. The Tables \ref{EXPON distribution Biases} and \ref{EXPON distribution MSE} below contain the results for the bias and MSE values.

\begin{table}[h!]
	\centering
	\setlength{\tabcolsep}{.4mm}
	\begin{tabular}{c||c||c|c|c||c|c|c|c|c}
		~~~$\vartheta_0$~~~ & ~~~$n$~~~ & ~~$\widehat{\vartheta}_n^{ML}$~~ & ~~$\widehat{\vartheta}_n^{MSE}$~~ & ~~$\widehat{\vartheta}_n^{CvM}$~~ & ~~$\widehat{\vartheta}_{n, 2}^{(0.25)}$~~ & ~~$\widehat{\vartheta}_{n, 2}^{(0.5)}$~~ & ~~$\widehat{\vartheta}_{n, 2}^{(1)}$~~ & ~~$\widehat{\vartheta}_{n, 2}^{(2)}$~~ & ~~$\widehat{\vartheta}_{n, 2}^{(3)}$~~ \\[1pt]
		\hline
		& 10 & 0.0557 & -0.0554 & 0.051 & 0.0428 & 0.0376 & 0.0333 & 0.0302 & \bf 0.0291 \\[1pt]
		\cline{2-10}
		& 25 & 0.0212 & -0.0205 & 0.0187 & 0.0161 & 0.0144 & 0.0129 & 0.0118 & \bf 0.0114 \\[1pt]
		\cline{2-10}
		0.5 & 50 & 0.0098 & -0.0106 & 0.0089 & 0.0075 & 0.0067 & 0.0062 & 0.0057 & \bf 0.0055 \\[1pt]
		\cline{2-10}
		& 100 & 0.005 & -0.0051 & 0.0045 & 0.0039 & 0.0035 & 0.0032 & 0.003 & \bf 0.0029 \\[1pt]
		\cline{2-10}
		& 200 & 0.0023 & -0.0027 & 0.002 & 0.0018 & 0.0016 & 0.0015 & 0.0014 & \bf 0.0013 \\[1pt]
		\hline
		& 10 & 0.2193 & -0.2245 & 0.2004 & 0.2011 & 0.1871 & 0.168 & 0.1476 & \bf 0.1371 \\[1pt]
		\cline{2-10}
		& 25 & 0.083 & -0.0836 & 0.0746 & 0.0754 & 0.0701 & 0.0634 & 0.0569 & \bf 0.0538 \\[1pt]
		\cline{2-10}
		2 & 50 & 0.0398 & -0.0418 & 0.0345 & 0.0358 & 0.0332 & 0.0298 & 0.0263 & \bf 0.0246 \\[1pt]
		\cline{2-10}
		& 100 & 0.0191 & -0.0213 & 0.0165 & 0.0171 & 0.0158 & 0.0142 & 0.0126 & \bf 0.0118 \\[1pt]
		\cline{2-10}
		& 200 & 0.0095 & -0.0106 & 0.0074 & 0.0084 & 0.0077 & 0.0067 & 0.0057 & \bf 0.0052 \\[1pt]
		\hline
		& 10 & 0.5437 & -0.5651 & 0.4863 & 0.5238 & 0.5059 & 0.4753 & 0.4303 & \bf 0.3997 \\[1pt]
		\cline{2-10}
		& 25 & 0.2102 & -0.2066 & 0.1832 & 0.2015 & 0.194 & 0.1818 & 0.1649 & \bf 0.154 \\[1pt]
		\cline{2-10}
		5 & 50 & 0.1048 & -0.0994 & 0.0923 & 0.1004 & 0.0967 & 0.0908 & 0.0829 & \bf 0.0779 \\[1pt]
		\cline{2-10}
		& 100 & 0.052 & -0.0491 & 0.044 & 0.0496 & 0.0477 & 0.0446 & 0.0404 & \bf 0.0378 \\[1pt]
		\cline{2-10}
		& 200 & 0.0264 & -0.0238 & 0.0224 & 0.0253 & 0.0243 & 0.0229 & 0.0209 & \bf 0.0196 \\[1pt]
		\hline
		& 10 & 1.123 & -1.1016 & 1.0316 & 1.1028 & 1.0837 & 1.0484 & 0.9885 & \bf 0.9401 \\[1pt]
		\cline{2-10}
		& 25 & 0.4177 & -0.4157 & 0.3719 & 0.4089 & 0.4008 & 0.3863 & 0.3628 & \bf 0.3448 \\[1pt]
		\cline{2-10}
		10 & 50 & 0.2041 & -0.204 & 0.1826 & 0.1996 & 0.1955 & 0.1883 & 0.1768 & \bf 0.1681 \\[1pt]
		\cline{2-10}
		& 100 & 0.0991 & -0.1029 & 0.0873 & 0.0967 & 0.0945 & 0.0908 & 0.0848 & \bf 0.0804 \\[1pt]
		\cline{2-10}
		& 200 & 0.0556 & \bf -0.045 & 0.0483 & 0.0544 & 0.0533 & 0.0513 & 0.0483 & 0.046 \\
	\end{tabular}
	\caption{Approximated biases calculated with 100,000 exponentially distributed Monte Carlo samples.} \label{EXPON distribution Biases}
\end{table}

\begin{table}[ht!]
	\centering
	\setlength{\tabcolsep}{.4mm}
	\begin{tabular}{c||c||c|c|c||c|c|c|c|c}
		~~~$\vartheta_0$~~~ & ~~~$n$~~~ & ~~$\widehat{\vartheta}_n^{ML}$~~ & ~~$\widehat{\vartheta}_n^{MSE}$~~ & ~~$\widehat{\vartheta}_n^{CvM}$~~ & ~~$\widehat{\vartheta}_n^{(0.25)}$~~ & ~~$\widehat{\vartheta}_n^{(0.5)}$~~ & ~~$\widehat{\vartheta}_n^{(1)}$~~ & ~~$\widehat{\vartheta}_n^{(2)}$~~ & ~~$\widehat{\vartheta}_n^{(3)}$~~ \\[1pt]
		\hline
		& 10 & 0.0416 & \bf 0.0277 & 0.0593 & 0.0409 & 0.0428 & 0.0496 & 0.0662 & 0.0837 \\[1pt]
		\cline{2-10}
		& 25 & 0.0123 & \bf 0.0105 & 0.0162 & 0.0127 & 0.0138 & 0.0167 & 0.0229 & 0.0294 \\[1pt]
		\cline{2-10}
		0.5 & 50 & 0.0055 & \bf 0.0051 & 0.0072 & 0.0058 & 0.0064 & 0.0078 & 0.0109 & 0.014 \\[1pt]
		\cline{2-10}
		& 100 & 0.0026 & \bf 0.0025 & 0.0034 & 0.0028 & 0.0031 & 0.0038 & 0.0053 & 0.0069 \\[1pt]
		\cline{2-10}
		& 200 & \bf 0.0013 & \bf 0.0013 & 0.0017 & 0.0014 & 0.0015 & 0.0019 & 0.0026 & 0.0034 \\[1pt]
		\hline
		& 10 & 0.6645 & \bf 0.4449 & 0.9504 & 0.6569 & 0.6525 & 0.6537 & 0.6845 & 0.7346 \\[1pt]
		\cline{2-10}
		& 25 & 0.1949 & \bf 0.1661 & 0.2573 & 0.1942 & 0.1952 & 0.2006 & 0.2184 & 0.2403 \\[1pt]
		\cline{2-10}
		2 & 50 & 0.0887 & \bf 0.0821 & 0.1165 & 0.0889 & 0.0898 & 0.0932 & 0.1029 & 0.1141 \\[1pt]
		\cline{2-10}
		& 100 & 0.0418 & \bf 0.0402 & 0.0549 & 0.042 & 0.0426 & 0.0444 & 0.0493 & 0.0549 \\[1pt]
		\cline{2-10}
		& 200 & 0.0205 & \bf 0.0201 & 0.0269 & 0.0207 & 0.021 & 0.022 & 0.0244 & 0.0272 \\[1pt]
		\hline
		& 10 & 4.0739 & \bf 2.7374 & 5.6848 & 4.0529 & 4.035 & 4.0092 & 3.9977 & 4.0335 \\[1pt]
		\cline{2-10}
		& 25 & 1.2302 & \bf 1.0465 & 1.621 & 1.2272 & 1.2259 & 1.2284 & 1.2493 & 1.2842 \\[1pt]
		\cline{2-10}
		5 & 50 & 0.5522 & \bf 0.5087 & 0.7246 & 0.5518 & 0.5524 & 0.5561 & 0.5706 & 0.5908 \\[1pt]
		\cline{2-10}
		& 100 & 0.2635 & \bf 0.2529 & 0.3445 & 0.2636 & 0.2641 & 0.2665 & 0.2745 & 0.285 \\[1pt]
		\cline{2-10}
		& 200 & 0.1295 & \bf 0.1268 & 0.1692 & 0.1296 & 0.1299 & 0.1313 & 0.1355 & 0.141 \\[1pt]
		\hline
		& 10 & 16.8106 & \bf 11.1652 & 23.5189 & 16.7647 & 16.7219 & 16.646 & 16.5344 & 16.4779 \\[1pt]
		\cline{2-10}
		& 25 & 4.885 & \bf 4.1598 & 6.4565 & 4.8785 & 4.8739 & 4.8702 & 4.8831 & 4.9188 \\[1pt]
		\cline{2-10}
		10 & 50 & 2.2069 & \bf 2.0371 & 2.8967 & 2.2053 & 2.2048 & 2.2069 & 2.2213 & 2.2464 \\[1pt]
		\cline{2-10}
		& 100 & 1.0473 & \bf 1.007 & 1.3747 & 1.0471 & 1.0474 & 1.0497 & 1.0594 & 1.074 \\[1pt]
		\cline{2-10}
		& 200 & 0.5126 & \bf 0.5014 & 0.6658 & 0.5127 & 0.513 & 0.5144 & 0.5198 & 0.5274 \\
	\end{tabular}
	\caption{Approximated MSE calculated with 100,000 exponentially distributed Monte Carlo samples.} \label{EXPON distribution MSE}
\end{table}

As for the biases, the maximum likelihood estimator and the minimum MSE estimator perform almost identically in terms of the absolute bias, and the minimum Cram\'{e}r-von Mises distance estimator has a slight edge. Our new estimator outperforms all other methods (virtually) uniformly. More precisely, it seems as if for larger tuning parameters '$a$' the bias decreases. We will show, however, that this observation is not correct in that generality. The results for the mean squared error reveal that the minimum MSE estimator is the best method with respect to this measure of quality, which is no surprise as it is constructed to minimize the MSE. For sample size $n = 10$ the superiority is particularly obvious, but for larger samples, the maximum likelihood estimator is only slightly worse. Our new estimator shows almost identical results (for $a = 0.25$) as the maximum likelihood estimator, undermining that the method is sound and powerful. In contrast to the observation with the bias values, the MSE appears to increase with '$a$'. This nicely illustrates the variance-bias trade-off commonly observed in the context of estimation problems.

\section{The case $\lowercase{a} \to \infty$}
\label{SEC the case a to infty}

As discussed previously, the simulation results for the exponential distribution somewhat indicate that as the tuning parameter '$a$' grows, the bias decreases while the MSE increases. Interestingly, we can lay observations for $a \to \infty$ on a rigorous theoretical basis. To be precise, observe the following general result.
\begin{theorem} \label{THM a to infinity}
	Consider the setting from Section \ref{SEC the new estimators} with weight function $w(t) = e^{-a t}$, $a > 0$. For the quantity $\psi_{n, q}(\vartheta, a) = \psi_{n, q}(\vartheta) = \lVert \eta_n(\, \cdot \, , \vartheta) \rVert_{L^q}$ from the end of Section \ref{SEC Existence and measurability}, we make the dependence on the tuning parameter '$a$' explicit. Then,
	\begin{align*}
	\lim_{a \, \to \, \infty} a^{q + 1} \, \big(\psi_{n, q}(\vartheta, a)\big)^q
	= \Gamma(q + 1) \left| \frac{1}{n} \sum_{j = 1}^{n} \frac{p^\prime_\vartheta(X_j)}{p_\vartheta(X_j)} \right|^q,
	\end{align*}
	on a set of measure one, where $\Gamma(\cdot)$ denotes the Gamma function.
\end{theorem}
The proof consists of an almost trivial application of an Abelian theorem for the Laplace transform, see p.182 of \cite{W:1959}, or the work by \cite{BGH:2000}. Since $a, q > 0$, the functions $\psi_{n, q}(\vartheta)$ and $a^{q + 1} \big( \psi_{n, q}(\vartheta) \big)^q$ attain their minimum in the same point. Thus, in the limit $a \to \infty$, our procedure essentially yields as an estimators the minimizer of the quantity
\begin{align*}
\left| \frac{1}{n} \sum_{j = 1}^{n} \frac{p^\prime_\vartheta(X_j)}{p_\vartheta(X_j)} \right|^q.
\end{align*}
In the situation of the exponential distribution as discussed in Section \ref{SEC exponential distribution}, the result reduces to $\lim_{a \, \to \, \infty} a^{3} \big( \psi_{n, 2}(\vartheta, a) \big)^2 = 2 \vartheta^2$, so in the limit $a \to \infty$, the procedure will choose $\widehat{\vartheta} = 0 \notin \Theta$ as the estimator, which leads to a bias of $- \vartheta_0$ and an MSE of $\vartheta_0^2$. The observation from the simulations is, therefore, not universal. An example for which the limit in Theorem \ref{THM a to infinity} is less trivial is the Rayleigh distribution.

\section{Example: Rayleigh distribution}
\label{SEC Rayleigh distribution}
Let $\Theta = (0, \infty)$ and take the density function of the Rayleigh distribution with parameter $\vartheta \in \Theta$,
\begin{align*}
p_\vartheta(x) = \frac{x}{\vartheta^2} \, \exp\Big( - \frac{x^2}{2 \vartheta^2} \Big), \quad x > 0.
\end{align*}
It is easy to check that the Rayleigh density satisfies all regularity conditions stated throughout the work, and that we have $\tfrac{p^\prime_\vartheta(x)}{p_\vartheta(x)} = \tfrac{1}{x} - \tfrac{x}{\vartheta^2}.$ The limit in Theorem \ref{THM a to infinity} thus takes the form
\begin{align*}
\Gamma(q + 1) \left| \frac{1}{n} \sum_{j = 1}^{n} \frac{p^\prime_\vartheta(X_j)}{p_\vartheta(X_j)} \right|^q
= \Gamma(q + 1) \left| \frac{1}{n} \sum_{j = 1}^{n} \left( \frac{1}{X_j} - \frac{X_j}{\vartheta^2} \right) \right|^q ,
\end{align*}
where $X_1, \dots, X_n$ are i.i.d. random variables which follow the Rayleigh law, $X_1 \sim p_{\vartheta_0}$, for some unknown scale parameter $\vartheta_0 \in \Theta$. In the case $q = 2$, it is easy to calculate that the minimum of the above function over $\vartheta > 0$ is given through
\begin{align*}
\widehat{\vartheta}_n^{AM} = \sqrt{\frac{\frac{1}{n} \sum_{j = 1}^{n} X_j}{\frac{1}{n} \sum_{j = 1}^{n} \frac{1}{X_j}}} .
\end{align*}
Strikingly, this asymptotically derived moment-type estimator is itself consistent for $\vartheta_0$, as
\begin{align*}
\widehat{\vartheta}_n^{AM} 
= \sqrt{\frac{\frac{1}{n} \sum_{j = 1}^{n} X_j}{\frac{1}{n} \sum_{j = 1}^{n} \frac{1}{X_j}}}
\longrightarrow \sqrt{\frac{\E \big[ X_1 \big]}{\E \Big[ \frac{1}{X_1} \Big]}}
= \sqrt{\frac{\sqrt{\frac{\pi}{2}} \, \vartheta_0}{\sqrt{\frac{\pi}{2}} \cdot \frac{1}{\vartheta_0}}}
= \vartheta_0
\end{align*}
$\mathbb{P}$-a.s., as $n \to \infty$, where we used the law of large numbers, as well as the fact that $X_1, \dots, X_n$ all follow the Rayleigh distribution with parameter $\vartheta_0$. We compare this estimator with other methods. Among them is our new estimator
\begin{align*}
\widehat{\vartheta}_{n, 2}^{(a)}
= \mathrm{arg \, min} \big\{ \big(\psi_{n, 2}(\vartheta)\big)^2 \, \big| \, \vartheta > 0 \big\} = \mathrm{arg \, min} \big\{ \vartheta^{-4} \widetilde{\Psi}_n^{(1)} + \vartheta^{-2} \widetilde{\Psi}_n^{(2)} + \widetilde{\Psi}_n^{(3)} \, \big| \, \vartheta > 0 \big\} ,
\end{align*}
where
\begin{align*}
\widetilde{\Psi}_n^{(1)}
&= \frac{2}{n^2} \sum_{1 \leq j < k \leq n} \left[ X_{(j)} X_{(k)} \cdot \frac{2}{a^3} \, \big( 1 - e^{- a X_{(j)}} \big) - \frac{X_{(j)}^2 X_{(k)}}{a^2} \, \big( e^{- a X_{(j)}} + e^{- a X_{(k)}} \big) \right] \\
&~~~ + \frac{1}{n^2} \sum_{j=1}^{n} \left[ \frac{2 X_{(j)}^2}{a^3} \, \big( 1 - e^{- a X_{(j)}} \big) - \frac{2 X_{(j)}^3}{a^2} \, e^{- a X_{(j)}} \right] ,
\end{align*}
\begin{align*}
\widetilde{\Psi}_n^{(2)}
&= \frac{2}{n^2} \sum_{1 \leq j < k \leq n} \left[ \frac{X_{(j)}^2 e^{- a X_{(k)}}}{a} \, \left( \frac{1}{a X_{(k)}} - 1 \right) + \frac{X_{(j)} e^{- a X_{(j)}}}{a} \, \left( \frac{X_{(j)}}{a X_{(k)}} - X_{(k)} \right) \right. \\
&\left.\phantom{\frac{X_{(j)}^2 e^{- a X_{(k)}}}{a}} ~~~~~~ - \frac{2}{a^3} \, \big( 1 - e^{- a X_{(j)}} \big) \left( \frac{X_{(k)}}{X_{(j)}} + \frac{X_{(j)}}{X_{(k)}} \right) \right] \\
&~~~ + \frac{1}{n^2} \sum_{j = 1}^{n} \left[ \frac{2 e^{-a X_{(j)}}}{a} \cdot X_{(j)} \left( \frac{2j}{a} - X_{(j)} \right) - \frac{4}{a^3} \, \big( 1 - e^{-a X_{(j)}} \big) \right] ,
\end{align*}
and
\begin{align*}
\widetilde{\Psi}_n^{(3)}
&= \frac{2}{n^2} \sum_{1 \leq j < k \leq n} \left[ \frac{X_{(j)}}{a X_{(k)}} \, e^{- a X_{(j)}} + \frac{2}{a^3 X_{(j)} X_{(k)}} \, \big( 1 - e^{- a X_{(j)}} \big) \right] \\
&~~~ + \frac{1}{n^2} \sum_{j=1}^{n} \left[ \frac{2}{a^3 X_{(j)}^2} \, \big( 1 - e^{- a X_{(j)}} \big) + \frac{e^{- a X_{(j)}}}{a} \, \left( 4j - 1 - \frac{2}{a X_{(j)}} \, (2j - 1) \right)  \right] ,
\end{align*}
and $X_{(1)} < \dotso < X_{(n)}$ denotes the ordered sample. It is easily seen that if both $\widetilde{\Psi}_n^{(1)} > 0$ and $\widetilde{\Psi}_n^{(2)} < 0$ $\mathbb{P}$-a.s., then the minimum can be calculated explicitly as
\begin{align*}
\widehat{\vartheta}_{n, 2}^{(a)} = \sqrt{- \frac{2 \widetilde{\Psi}_n^{(1)}}{\widetilde{\Psi}_n^{(2)}}},
\end{align*}
and indeed, using that $e^{- a X_{(k)}} < e^{- a X_{(j)}}$  and $1 - e^{- a X_{(j)}} - a X_{(j)} e^{- a X_{(j)}} > 0$ $\mathbb{P}$-a.s., we have
\begin{align*}
\widetilde{\Psi}_n^{(1)}
&> \frac{2}{n^2} \sum_{1 \leq j < k \leq n} \frac{2 X_{(j)} X_{(k)}}{a^3} \Big( 1 - e^{- a X_{(j)}} - a X_{(j)} e^{- a X_{(j)}} \Big) + \frac{1}{n^2} \sum_{j = 1}^{n} \frac{2 X_{(j)}^2}{a^3} \Big( 1 - e^{- a X_{(j)}} - a X_{(j)} e^{- a X_{(j)}} \Big) \\
&> 0 \quad \mathbb{P}-\text{a.s.} ,
\end{align*}
and with similar thoughts, $\widetilde{\Psi}_n^{(2)} < 0$ $\mathbb{P}$-a.s.. Additionally, we consider the maximum likelihood estimator and a moment estimator, which are given as
\begin{align*}
\widehat{\vartheta}_n^{ML} = \sqrt{\frac{1}{2n} \sum_{j = 1}^{n} X_j^2} \quad \text{and} \quad
\widehat{\vartheta}_n^{Mom} = \sqrt{\frac{2}{\pi}} \cdot \frac{1}{n} \sum_{j = 1}^{n} X_j ,
\end{align*}
respectively. Note in particular that the moment estimator is unbiased and we can expect it to outperform the other estimators in this regard. Finally, we include the minimum Cram\'{e}r-von Mises distance estimator given through
\begin{align*}
\widehat{\vartheta}_n^{CvM}
= \arg\min \left\{ \frac{1}{n} \sum_{j = 1}^{n} \left[ \left( \frac{2j - 1}{n} - 2 \right) \exp\left( - \frac{X_{(j)}^2}{2 \vartheta^2} \right) + \exp\left( - \frac{X_{(j)}^2}{\vartheta^2} \right) \right] ~ \Big| ~ \vartheta > 0 \right\} ,
\end{align*}
where we solve the minimization numerically via a sequential least squares programming method as in the case of the exponential distribution in Section \ref{SEC exponential distribution}, using as initial value the maximum likelihood estimator. The execution of the comparison is as in the example on the exponential distribution, and the results are displayed in Tables \ref{RAYLEIGH distribution Biases} and \ref{RAYLEIGH distribution MSE}.

\begin{table}[h!]
	\centering
	\setlength{\tabcolsep}{.4mm}
	\begin{tabular}{c||c||c|c|c|c||c|c|c|c|c}
		~~~$\vartheta_0$~~~ & ~~~$n$~~~ & ~~$\widehat{\vartheta}_n^{ML}$~~ & ~~$\widehat{\vartheta}_n^{Mom}$~~ & ~~$\widehat{\vartheta}_n^{AM}$~~ & ~~$\widehat{\vartheta}_n^{CvM}$~~ & ~~$\widehat{\vartheta}_{n, 2}^{(0.25)}$~~ & ~~$\widehat{\vartheta}_{n, 2}^{(0.5)}$~~ & ~~$\widehat{\vartheta}_{n, 2}^{(1)}$~~ & ~~$\widehat{\vartheta}_{n, 2}^{(2)}$~~ & ~~$\widehat{\vartheta}_{n, 2}^{(3)}$~~ \\[1pt]
		\hline
		& 10 & -0.0061 & \bf 0.0001 & 0.0269 & 0.0029 & -0.005 & -0.004 & -0.0023 & 0.0006 & 0.0028 \\[1pt]
		\cline{2-11}
		& 25 & -0.0025 & \bf 0.0001 & 0.014 & 0.0012 & -0.002 & -0.0016 & -0.0009 & 0.0002 & 0.0011 \\[1pt]
		\cline{2-11}
		0.5 & 50 & -0.0011 & \bf 0.0001 & 0.0085 & 0.0007 & -0.0009 & -0.0007 & -0.0003 & 0.0002 & 0.0007 \\[1pt]
		\cline{2-11}
		& 100 & -0.0006 & \bf 0 & 0.0047 & 0.0003 & -0.0005 & -0.0004 & -0.0002 & 0.0001  & 0.0003 \\[1pt]
		\cline{2-11}
		& 200 & -0.0004 & \bf 0 & 0.0027 & 0.0001 & -0.0003 & -0.0002 & -0.0002 & \bf 0 & 0.0001 \\[1pt]
		\hline
		& 10 & -0.0246 & \bf 0.0001 & 0.1074 & 0.0107 & -0.0093 & 0.0019 & 0.0189 & 0.0434 & 0.0598 \\[1pt]
		\cline{2-11}
		& 25 & -0.0097 & \bf 0 & 0.0549 & 0.0039 & -0.0036 & 0.0007 & 0.0073 & 0.0169 & 0.0237 \\[1pt]
		\cline{2-11}
		2 & 50 & -0.0057 & -0.0007 & 0.032 & 0.0013 & -0.0025 & \bf -0.0003 & 0.003 & 0.0079 & 0.0114 \\[1pt]
		\cline{2-11}
		& 100 & -0.0026 & -0.0001 & 0.019 & 0.0008 & -0.0011 & \bf 0 & 0.0017 & 0.0041 & 0.0058 \\[1pt]
		\cline{2-11}
		& 200 & -0.0015 & -0.0003 & 0.0104 & \bf 0 & -0.0007 & -0.0002 & 0.0005 & 0.0016 & 0.0025 \\[1pt]
		\hline
		& 10 & -0.0624 & \bf -0.0009 & 0.2642 & 0.0255 & 0.0156 & 0.064 & 0.1293 & 0.1926 & 0.2199 \\[1pt]
		\cline{2-11}
		& 25 & -0.0245 & \bf -0.0002 & 0.1388 & 0.0097 & 0.0063 & 0.0251 & 0.0519 & 0.0817 & 0.0973 \\[1pt]
		\cline{2-11}
		5 & 50 & -0.0132 & \bf -0.0002 & 0.0848 & 0.0049 & 0.003 & 0.0129 & 0.027 & 0.0432 & 0.0523 \\[1pt]
		\cline{2-11}
		& 100 & -0.0059 & \bf 0 & 0.0477 & 0.0021 & 0.0016 & 0.0062 & 0.0128 & 0.0206 & 0.0253 \\[1pt]
		\cline{2-11}
		& 200 & -0.0028 & \bf 0.0002 & 0.0279 & 0.001 & 0.001 & 0.0033 & 0.0066 & 0.0106 & 0.0129 \\[1pt]
		\hline
		& 10 & -0.1248 & \bf -0.0004 & 0.5383 & 0.0537 & 0.1302 & 0.2617 & 0.3919 & 0.4777 & 0.5055 \\[1pt]
		\cline{2-11}
		& 25 & -0.0565 & \bf -0.0076 & 0.2699 & 0.0123 & 0.043 & 0.0965 & 0.1564 & 0.2074 & 0.2293 \\[1pt]
		\cline{2-11}
		10 & 50 & -0.0261 & \bf -0.0021 & 0.1582 & 0.0083 & 0.0225 & 0.048 & 0.0783 & 0.1073 & 0.1214 \\[1pt]
		\cline{2-11}
		& 100 & -0.0109 & \bf 0.0013 & 0.0979 & 0.0057 & 0.0138 & 0.0272 & 0.043 & 0.0586 & 0.0671 \\[1pt]
		\cline{2-11}
		& 200 & -0.0077 & \bf -0.001 & 0.0545 & 0.0011 & 0.0057 & 0.0128 & 0.0207 & 0.0289 & 0.0334 \\
	\end{tabular}
	\caption{Approximated biases calculated with 100,000 Rayleigh-distributed Monte Carlo samples.} \label{RAYLEIGH distribution Biases}
\end{table}

\begin{table}[ht!]
	\centering
	\setlength{\tabcolsep}{.4mm}
	\begin{tabular}{c||c||c|c|c|c||c|c|c|c|c}
		~~~$\vartheta_0$~~~ & ~~~$n$~~~ & ~~$\widehat{\vartheta}_n^{ML}$~~ & ~~$\widehat{\vartheta}_n^{Mom}$~~ & ~~$\widehat{\vartheta}_n^{AM}$~~ & ~~$\widehat{\vartheta}_n^{CvM}$~~ & ~~$\widehat{\vartheta}_{n, 2}^{(0.25)}$~~ & ~~$\widehat{\vartheta}_{n, 2}^{(0.5)}$~~ & ~~$\widehat{\vartheta}_{n, 2}^{(1)}$~~ & ~~$\widehat{\vartheta}_{n, 2}^{(2)}$~~ & ~~$\widehat{\vartheta}_{n, 2}^{(3)}$~~ \\[1pt]
		\hline
		& 10 & \bf 0.0061 & 0.0067 & 0.0135 & 0.0082 & 0.0062 & 0.0062 & 0.0064 & 0.0068 & 0.0072 \\[1pt]
		\cline{2-11}
		& 25 & \bf 0.0025 & 0.0027 & 0.0061 & 0.0033 & \bf 0.0025 & \bf 0.0025 & \bf 0.0025 & 0.0027 & 0.0028 \\[1pt]
		\cline{2-11}
		0.5 & 50 & \bf 0.0013 & 0.0014 & 0.0033 & 0.0017 & \bf 0.0013 & \bf 0.0013 & \bf 0.0013 & \bf 0.0013 & 0.0014 \\[1pt]
		\cline{2-11}
		& 100 & \bf 0.0006 & 0.0007 & 0.0019 & 0.0008 & \bf 0.0006 & \bf 0.0006 & \bf 0.0006 & 0.0007 & 0.0007 \\[1pt]
		\cline{2-11}
		& 200 & \bf 0.0003 & \bf 0.0003 & 0.001 & 0.0004 & \bf 0.0003 & \bf 0.0003 & \bf 0.0003 & \bf 0.0003 & 0.0004 \\[1pt]
		\hline
		& 10 & \bf 0.0992 & 0.1089 & 0.2185 & 0.1316 & 0.103 & 0.109 & 0.1248 & 0.157 & 0.18 \\[1pt]
		\cline{2-11}
		& 25 & \bf 0.0401 & 0.0438 & 0.0977 & 0.0527 & 0.0412 & 0.0433 & 0.0488 & 0.0601 & 0.069 \\[1pt]
		\cline{2-11}
		2 & 50 & \bf 0.0199 & 0.0219 & 0.0543 & 0.0264 & 0.0205 & 0.0215 & 0.0242 & 0.0297 & 0.0341 \\[1pt]
		\cline{2-11}
		& 100 & \bf 0.01 & 0.0109 & 0.03 & 0.013 & 0.0102 & 0.0106 & 0.0119 & 0.0146 & 0.0167 \\[1pt]
		\cline{2-11}
		& 200 & \bf 0.005 & 0.0055 & 0.0167 & 0.0065 & 0.0051 & 0.0053 & 0.006 & 0.0073 & 0.0083 \\[1pt]
		\hline
		& 10 & \bf 0.6205 & 0.6827 & 1.3695 & 0.8271 & 0.7057 & 0.8359 & 1.0635 & 1.2775 & 1.3473 \\[1pt]
		\cline{2-11}
		& 25 & \bf 0.2521 & 0.276 & 0.6122 & 0.3314 & 0.2803 & 0.3258 & 0.4097 & 0.5088 & 0.5566 \\[1pt]
		\cline{2-11}
		5 & 50 & \bf 0.125 & 0.1371 & 0.3398 & 0.1648 & 0.1385 & 0.1606 & 0.2015 & 0.2529 & 0.2811 \\[1pt]
		\cline{2-11}
		& 100 & \bf 0.0627 & 0.0684 & 0.1876 & 0.0819 & 0.0688 & 0.0793 & 0.0989 & 0.1242 & 0.1392 \\[1pt]
		\cline{2-11}
		& 200 & \bf 0.0311 & 0.0341 & 0.1039 & 0.0409 & 0.0343 & 0.0395 & 0.0491 & 0.0618 & 0.0695 \\[1pt]
		\hline
		& 10 & \bf 2.4749 & 2.7278 & 5.4528 & 3.3202 & 3.3485 & 4.2534 & 5.0993 & 5.4722 & 5.523 \\[1pt]
		\cline{2-11}
		& 25 & \bf 0.9966 & 1.0933 & 2.4305 & 1.3171 & 1.2922 & 1.6219 & 2.0109 & 2.3024 & 2.3989 \\[1pt]
		\cline{2-11}
		10 & 50 & \bf 0.5 & 0.5445 & 1.3528 & 0.6504 & 0.6342 & 0.7926 & 0.9955 & 1.1792 & 1.26 \\[1pt]
		\cline{2-11}
		& 100 & \bf 0.2499 & 0.2735 & 0.7504 & 0.329 & 0.3179 & 0.3959 & 0.4966 & 0.5962 & 0.6473 \\[1pt]
		\cline{2-11}
		& 200 & \bf 0.1248 & 0.1364 & 0.42 & 0.1637 & 0.1579 & 0.1961 & 0.247 & 0.2999 & 0.3288 \\
	\end{tabular}
	\caption{Approximated MSE calculated with 100,000 Rayleigh-distributed Monte Carlo samples.} \label{RAYLEIGH distribution MSE}
\end{table}

\noindent
Apparently, the moment estimator $\widehat{\vartheta}_n^{Mom}$ outperforms the other estimators with respect to the bias values, while the maximum likelihood estimator $\widehat{\vartheta}_n^{ML}$ gets the smallest MSE. The estimator we obtained via the limit results from the previous section seems sound in itself but is completely negligible compared to the other methods. In terms of bias, the minimum Cram\'{e}r-von Mises distance estimator is preferable to the maximum likelihood method, and both are outdone by our new estimator, which even keeps up with the unbiased moment estimator for the smaller values of the parameter $\vartheta_0$. Notice that the maximum likelihood and moment estimator tend to underestimate the parameter, while the other procedures tend to a slight overestimation. As for the MSE, the moment estimator and our new method perform similarly and follow the maximum likelihood estimator closely. The minimum Cram\'{e}r-von Mises distance estimator is a bit behind. To summarize, the maximum likelihood and moment estimator for the Rayleigh parameter are both simple and very convincing, but the newly proposed method keeps up (for suitably chosen tuning parameter) and appears to find a good compromise between bias and MSE. The only graver weakness shows for the large parameter value $\vartheta_0 = 10$ and small sample sizes $n = 10, 25$.

\section{Example: The Burr Type XII distribution}
\label{SEC Burr distribution}

Consider the density function $p_\vartheta(x) = c \, k \, x^{c - 1} \, \big( 1 + x^c \big)^{- k - 1}$, $x > 0$, where $\vartheta = (c, k) \in (0, \infty)^2 = \Theta$. It is not exactly trivial, but still straight-forward, to prove that this is an admissible distribution in terms of the setting in Section \ref{SEC the new estimators} [see also \cite{BE:2018}] and the conditions (R1) -- (R3). With $q = 2$ and weight $w(t) = e^{- a t}$, where $a > 0$, the function $\psi_{n, 2}(\vartheta) = \big\lVert \eta_n( \, \cdot \, , \vartheta) \big\rVert_{L^2}$ from Section \ref{SEC Existence and measurability} (see also Section \ref{SEC the new estimators}) can be calculated explicitly as
\begin{align*}
\Big(\psi_{n, 2}(\vartheta)\Big)^2
&= \frac{2}{n^2} \sum_{1 \leq j < \ell \leq n} \left\{ A_{(\ell)}(c, k) \left[ \frac{2 A_{(j)}(c, k)}{a^3} \, \big( 1 - e^{- a X_{(j)}} \big) + \frac{B_{(j)}(c, k)}{a^2} \, \big( e^{- a X_{(j)}} + e^{- a X_{(\ell)}} \big) \phantom{\frac{X_{(j)}}{a}} \right.\right. \\
&\left. \left. \qquad\qquad\qquad\qquad\qquad\quad~~+ \frac{c - 2}{a^2} \, e^{- a X_{(j)}} - \frac{X_{(j)}}{a} \, e^{- a X_{(j)}} \right] + \frac{B_{(j)}(c, k)}{a} \, e^{- a X_{(\ell)}} \right\} \\
&~~~ + \frac{1}{n^2} \sum_{j = 1}^{n} \left\{ \big(A_{(j)}(c, k)\big)^2 \left( - \frac{2 X_{(j)}}{a^2} \, e^{- a X_{(j)}} - \frac{2}{a^3} \, e^{- a X_{(j)}} + \frac{2}{a^3} \right) \right. \\
&\left. \phantom{\frac{2 X_{(j)}}{a^2}}  ~~~~~~~~~~~+ \frac{2(j - 1) \, c}{a^2} \, A_{(j)}(c, k) \, e^{- a X_{(j)}} + \frac{2 B_{(j)}(c, k)}{a} \, e^{- a X_{(j)}} \right\} \\
&~~~ + \frac{2 c}{a \, n^2} \sum_{j = 1}^{n} j \, e^{- a X_{(j)}} - \frac{1}{a \, n^2} \sum_{j = 1}^{n} e^{- a X_{(j)}} ,
\end{align*}
where 
\begin{align*}
A_{(j)}(c, k) = c \, (k + 1) \, \frac{X_{(j)}^{c - 1}}{1 + X_{(j)}^c} - \frac{c - 1}{X_{(j)}}, \quad
B_{(j)}(c, k) = - c \, (k + 1) \, \frac{X_{(j)}^c}{1 + X_{(j)}^c} ,
\end{align*}
and where $X_{(1)} < \dotso < X_{(n)}$ denotes the ordered sample. Our estimator $\widehat{\vartheta}_{n, 2}^{(a)} = \big( \widehat{c}_n^{(a)}, \, \widehat{k}_n^{(a)} \big)$, as defined in (\ref{Estimator of theta}), can be calculated as the minimizer of the above function over $\Theta$. We use the 'L-BFGS-B'-method [L-BFGS-B algorithm, see \cite{BLNZ:1995} and \cite{ZBLN:1997}] implemented in the 'optimize.minimize' function of 'scipy' to solve the minimization numerically, using $(1, 1)$ as initial values. (Note that in preliminary simulations we have tried several other optimization routines, like a truncated Newton algorithm or the 'SLSQP' from previous sections, but the 'L-BFGS-B'-method appeared to be the most reliable for our purpose.) 

\begin{table}[h!]
	\centering
	\setlength{\tabcolsep}{.4mm}
	\begin{tabular}{c||c||c|c||c|c|c|c|c}
		~~~$\vartheta_0 = {c_0 \choose k_0}$~~~ & ~~~$n$~~~ & ~~$\widehat{\vartheta}_n^{ML}$~~ & ~~$\widehat{\vartheta}_n^{CvM}$~~ & ~~$\widehat{\vartheta}_{n, 2}^{(0.25)}$~~ & ~~$\widehat{\vartheta}_{n, 2}^{(0.5)}$~~ & ~~$\widehat{\vartheta}_{n, 2}^{(1)}$~~ & ~~$\widehat{\vartheta}_{n, 2}^{(2)}$~~ & ~~$\widehat{\vartheta}_{n, 2}^{(3)}$~~ \\
		\hline
		& 10 & -- & 0.142 & \bf 0.0094 & -0.1608 & 0.0375 & 0.0636 & 0.0585 \\[1pt]
		\cline{3-9}
		& & -- & 1.3014 & 0.175 & -0.2579 & 0.0243 & 0.1635 & \bf 0.1382 \\[1pt]
		\cline{2-9}
		& 25 & 0.0406 & 0.0451 & -0.0377 & -0.196 & 0.024 & 0.0223 & \bf 0.0152 \\[1pt]
		\cline{3-9}
		& & 0.1057 & 0.1311 & 0.0184 & -0.301 & 0.0503 & 0.049 & \bf 0.0131 \\[1pt]
		\cline{2-9}
		${\, 0.8 \, \choose 2}$ & 50 & 0.0197 & 0.0207 & -0.0411 & -0.1673 & 0.0121 & 0.0102 & \bf 0.0034 \\[1pt]
		\cline{3-9}
		& & 0.0491 & 0.0565 & \bf -0.0089 & -0.2435 & 0.0256 & 0.0214 & -0.0122 \\[1pt]
		\cline{2-9}
		& 100 & 0.0097 & 0.0102 & -0.0307 & -0.1102 & 0.0062 & 0.0056 & \bf 0.001 \\[1pt]
		\cline{3-9}
		& & 0.0234 & 0.0266 & -0.0128 & -0.1505 & 0.0125 & 0.0125 & \bf -0.0095 \\[1pt]
		\cline{2-9}
		& 200 & 0.0046 & 0.005 & -0.012 & -0.051 & 0.003 & 0.0029 & \bf 0.0012 \\[1pt]
		\cline{3-9}
		& & 0.0114 & 0.0131 & -0.0048 & -0.066 & 0.0064 & 0.0071 & \bf -0.0001 \\[1pt]
		\hline
		& 10 & 0.2956 & 0.3458 & -0.2755 & -1.0773 & 0.2152 & 0.1987 & \bf 0.1841 \\[1pt]
		\cline{3-9}
		& & 2.8551 & 37.1075 & \bf 1.6188 & -2.1985 & 2.35 & 2.287 & 2.1681 \\[1pt]
		\cline{2-9}
		& 25 & 0.1027 & 0.1082 & -0.1434 & -1.2772 & 0.0725 & 0.0655 & \bf 0.0618 \\[1pt]
		\cline{3-9}
		& & 0.6208 & 0.8619 & \bf 0.2341 & -2.6011 & 0.5033 & 0.4754 & 0.4647 \\[1pt]
		\cline{2-9}
		${\, 2 \, \choose 5}$ & 50 & 0.0476 & 0.0492 & -0.0347 & -1.4268 & 0.0326 & 0.0298 & \bf 0.0283 \\[1pt]
		\cline{3-9}
		& & 0.2669 & 0.3415 & \bf 0.1278 & -2.809 & 0.2126 & 0.2039 & 0.2021 \\[1pt]
		\cline{2-9}
		& 100 & 0.0233 & 0.0233 & \bf 0.0079 & -1.5877 & 0.0159 & 0.0145 & 0.0138 \\[1pt]
		\cline{3-9}
		& & 0.1285 & 0.1565 & \bf 0.0946 & -3.0394 & 0.1025 & 0.0989 & 0.0983 \\[1pt]
		\cline{2-9}
		& 200 & 0.012 & 0.0113 & 0.0089 & -1.7622 & 0.0082 & 0.0076 & \bf 0.0073 \\[1pt]
		\cline{3-9}
		& & 0.0627 & 0.0732 & 0.0526 & -3.3064 & 0.05 & 0.0485 & \bf 0.0483 \\[1pt]
		\hline
		& 10 & -- & 2.1411 & \bf 1.0267 & 1.0622 & 1.0927 & 1.1037 & 1.1374 \\[1pt]
		\cline{3-9}
		& & -- & 0.0451 & 0.0322 & \bf 0.0313 & 0.0327 & 0.038 & 0.0432 \\[1pt]
		\cline{2-9}
		& 25 & 0.3731 & 0.4635 & 0.3177 & 0.3233 & 0.3167 & \bf 0.3101 & 0.3146 \\[1pt]
		\cline{3-9}
		& & 0.0143 & 0.0113 & 0.0096 & \bf 0.0095 & 0.0106 & 0.013 & 0.0151 \\[1pt]
		\cline{2-9}
		${5 \choose \, 0.8 \,}$ & 50 & 0.1748 & 0.2071 & 0.1519 & 0.153 & 0.1488 & \bf 0.1453 & 0.1463 \\[1pt]
		\cline{3-9}
		& & 0.0063 & 0.0046 & \bf 0.0038 & 0.0039 & 0.0045 & 0.0057 & 0.0067 \\[1pt]
		\cline{2-9}
		& 100 & 0.0835 & 0.096 & 0.0731 & 0.0729 & 0.0708 & \bf 0.069 & 0.0693 \\[1pt]
		\cline{3-9}
		& & 0.0031 & 0.0023 & \bf 0.0019 & \bf 0.0019 & 0.0022 & 0.0028 & 0.0033 \\[1pt]
		\cline{2-9}
		& 200 & 0.0421 & 0.0481 & 0.0375 & 0.037 & 0.036 & 0.0348 & \bf 0.0347 \\[1pt]
		\cline{3-9}
		& & 0.0016 & 0.0012 & \bf 0.001 & \bf 0.001 & 0.0012 & 0.0015 & 0.0017 \\[1pt]
	\end{tabular}
	\caption{Approximated biases calculated with 100,000 Burr-distributed Monte Carlo samples.} \label{BURR distribution Biases}
\end{table}

As competitors to our estimator we consider the maximum likelihood estimator with implementation as suggested by \cite{SG:1993} [for a different algorithm, see \cite{W:1983}]. More precisely we use the Newton-Raphson method (with initial value $c = 1$) to find the root
\begin{align*}
\frac{n}{c} + \sum_{j = 1}^{n} \log(X_j) - \left[ \left( \frac{1}{n} \sum_{j = 1}^{n} \log(1 + X_j^c) \right)^{-1} + 1 \right] \cdot \sum_{j = 1}^{n} \frac{X_j^c \, \log(X_j)}{1 + X_j^c} \stackrel{!}{=} 0
\end{align*}
giving an estimate $\widehat{c}_n^{ML}$ for $c$ which we then introduce into
\begin{align*}
\widehat{k}_n^{ML} = \left( \frac{1}{n} \sum_{j = 1}^{n} \log(1 + X_j^{\widehat{c}_n^{ML}}) \right)^{-1} .
\end{align*}
Both relations are easily derived from the likelihood equations. Additionally, we consider the minimum Cram\'{e}r-von Mises distance estimator, which can be calculated from
\begin{align*}
\widehat{\vartheta}_n^{CvM} 
= \big( \widehat{c}_n^{CvM}, \, \widehat{k}_n^{CvM} \big)
= \arg\min \left\{ \frac{1}{n} \sum_{j = 1}^{n} \big( 1 + X_{(j)}^c \big)^{- k} \left[ \frac{2 j - 1}{n} - 2 + \big( 1 + X_{(j)}^c \big)^{- k} \right] ~ \Big| ~ c, k > 0 \right\}
\end{align*}
(the minimization is solved numerically, similar to our new estimator). Note that there have been further contributions to the estimation of the Burr parameters [see \cite{S:1983}, \cite{SG:1993}, \cite{W:1993}, and \cite{WC:2010}]. 

\begin{table}[hb!]
	\centering
	\setlength{\tabcolsep}{.4mm}
	\begin{tabular}{c||c||c|c||c|c|c|c|c}
		~~~$\vartheta_0 = {c_0 \choose k_0}$~~~ & ~~~$n$~~~ & ~~$\widehat{\vartheta}_n^{ML}$~~ & ~~$\widehat{\vartheta}_n^{CvM}$~~ & ~~$\widehat{\vartheta}_{n, 2}^{(0.25)}$~~ & ~~$\widehat{\vartheta}_{n, 2}^{(0.5)}$~~ & ~~$\widehat{\vartheta}_{n, 2}^{(1)}$~~ & ~~$\widehat{\vartheta}_{n, 2}^{(2)}$~~ & ~~$\widehat{\vartheta}_{n, 2}^{(3)}$~~ \\
		\hline
		& 10 & -- & 0.179 & 0.1398 & 0.2562 & 0.1209 & \bf 0.1092 & 0.1199 \\[1pt]
		\cline{3-9}
		& & -- & 22686.9274 & 1.9083 & \bf 0.9712 & 1.0083 & 2.2673 & 2.5852 \\[1pt]
		\cline{2-9}
		& 25 & \bf 0.0228 & 0.0293 & 0.0724 & 0.2017 & 0.031 & 0.0344 & 0.0406 \\[1pt]
		\cline{3-9}
		& & \bf 0.2354 & 0.3578 & 0.2867 & 0.3025 & 0.3438 & 0.513 & 0.7007 \\[1pt]
		\cline{2-9}
		${\, 0.8 \, \choose 2}$ & 50 & \bf 0.009 & 0.0121 & 0.0518 & 0.1577 & 0.0137 & 0.0159 & 0.02 \\[1pt]
		\cline{3-9}
		& & \bf 0.0957 & 0.1242 & 0.1437 & 0.1853 & 0.1557 & 0.2298 & 0.343 \\[1pt]
		\cline{2-9}
		& 100 & \bf 0.0042 & 0.0056 & 0.0337 & 0.1009 & 0.0066 & 0.0076 & 0.01 \\[1pt]
		\cline{3-9}
		& & \bf 0.0483 & 0.0545 & 0.0781 & 0.1136 & 0.0737 & 0.1066 & 0.171 \\[1pt]
		\cline{2-9}
		& 200 & \bf 0.002 & 0.0027 & 0.0144 & 0.0472 & 0.0032 & 0.0036 & 0.0045 \\[1pt]
		\cline{3-9}
		& & 0.207 & \bf 0.0255 & 0.0368 & 0.0582 & 0.0358 & 0.0507 & 0.0769 \\[1pt]
		\hline
		& 10 & 0.4819 & 0.9394 & 1.0985 & 2.6895 & 0.4657 & 0.4352 & \bf 0.4291 \\[1pt]
		\cline{3-9}
		& & 260.2839 & 1671094.1072 & 88.3851 & \bf 22.0669 & 164.7637 & 180.7299 & 168.6383 \\[1pt]
		\cline{2-9}
		& 25 & 0.1139 & 0.1717 & 0.4778 & 2.7278 & 0.1143 & \bf 0.1135 & 0.1163 \\[1pt]
		\cline{3-9}
		& & 3.5036 & 10.6678 & 4.6982 & 9.531 & 3.3494 & \bf 3.3022 & 3.3652 \\[1pt]
		\cline{2-9}
		${\, 2 \, \choose 5}$ & 50 & \bf 0.0477 & 0.0699 & 0.1719 & 2.9353 & 0.0497 & 0.0504 & 0.0526 \\[1pt]
		\cline{3-9}
		& & 1.059 & 1.8952 & 1.5019 & 10.0855 & \bf 1.0549 & 1.089 & 1.1556 \\[1pt]
		\cline{2-9}
		& 100 & \bf 0.0221 & 0.0321 & 0.041 & 3.216 & 0.0237 & 0.0243 & 0.0254 \\[1pt]
		\cline{3-9}
		& & \bf 0.4312 & 0.694 & 0.5032 & 10.9842 & 0.4435 & 0.4665 & 0.5016 \\[1pt]
		\cline{2-9}
		& 200 & \bf 0.0107 & 0.0153 & 0.0119 & 3.5414 & 0.0116 & 0.0119 & 0.0125 \\[1pt]
		\cline{3-9}
		& & \bf 0.1954 & 0.2996 & 0.2009 & 12.0603 & 0.2039 & 0.2156 & 0.233 \\[1pt]
		\hline
		& 10 & -- & 63.9726 & \bf 14.0301 & 14.8418 & 16.3923 & 16.7479 & 17.2525 \\[1pt]
		\cline{3-9}
		& & -- & 0.2941 & 0.1349 & \bf 0.137 & 0.1405 & 0.1479 & 0.1565 \\[1pt]
		\cline{2-9}
		& 25 & 1.6759 & 2.7939 & 1.6869 & 1.6671 & 1.631 & \bf 1.6163 & 1.6767 \\[1pt]
		\cline{3-9}
		& & \bf 0.0403 & 0.0496 & 0.0443 & 0.045 & 0.0461 & 0.0483 & 0.0508 \\[1pt]
		\cline{2-9}
		${5 \choose \, 0.8 \,}$ & 50 & \bf 0.6177 & 0.8648 & 0.6609 & 0.646 & 0.6299 & 0.6307 & 0.6597 \\[1pt]
		\cline{3-9}
		& & \bf 0.0189 & 0.0227 & 0.0212 & 0.0215 & 0.022 & 0.0229 & 0.0239 \\[1pt]
		\cline{2-9}
		& 100 & \bf 0.2746 & 0.3641 & 0.2958 & 0.2891 & 0.2827 & 0.2844 & 0.2992 \\[1pt]
		\cline{3-9}
		& & \bf 0.0091 & 0.0108 & 0.0103 & 0.0105 & 0.0107 & 0.0111 & 0.0116 \\[1pt]
		\cline{2-9}
		& 200 & \bf 0.1293 & 0.1671 & 0.1395 & 0.1366 & 0.1338 & 0.135 & 0.143 \\[1pt]
		\cline{3-9}
		& & \bf 0.0045 & 0.0054 & 0.0052 & 0.0052 & 0.0053 & 0.0055 & 0.0058 \\[1pt]
	\end{tabular}
	\caption{Approximated MSE values calculated with 100,000 Burr-distributed Monte Carlo samples.} \label{BURR distribution MSE}
\end{table}

Like for the exponential- and Rayleigh distribution, we approximate bias and MSE of these estimators and show the results in Tables \ref{BURR distribution Biases} and \ref{BURR distribution MSE}. For each value of $\vartheta_0$ and $n$, the first line corresponds to the bias/MSE of the estimator for the $c$-parameter, and the second line corresponds to the $k$-parameter. As before, it becomes evident that our new procedure outperforms the maximum likelihood and minimum Cram\'{e}r-von Mises distance estimator in terms of the bias. Unlike for the exponential distribution, the dependence on the tuning parameter '$a$' is less clear: For a great deal of parameter values and sample sizes, the estimator $\widehat{\vartheta}_{n, 2}^{(3)}$ yields the best result, but in some cases (mostly for the $k$-parameter) the estimator $\widehat{\vartheta}_{n, 2}^{(0.25)}$, with tuning parameter from the other end of the spectrum, performs best. Also observe the oddity that in some cases the estimator fares noticeably worse for $a = 0.5$ than for both smaller and larger tuning parameters. Thus, if one seeks to minimize some measure of quality of the estimators, an optimal, data dependent choice of the tuning parameter would be useful (more on this in Section \ref{SEC notes and comments}). In the light of our simulations, we suggest the use of $\widehat{\vartheta}_{n, 2}^{(3)}$ in practice as long as no adaptive tuning is available. Both in the bias and in the MSE simulation, the maximum likelihood estimator ran into computational issues for sample size $n = 10$. The minimum Cram\'{e}r-von Mises distance estimator is more stable in this regard, but still a lot less so than our new estimators which show notably slighter outliers only for large values of the Burr parameters. Once samples get larger ($n = 50+$), the asymptotic optimality properties of the maximum likelihood estimator appear to kick in, as its performance stabilizes. Still for suitably chosen tuning parameter, our estimators are very close in virtually all instances. The small sample behavior of the maximum likelihood estimator poses a huge drawback for applications and the problem is well-known.

\section{Example: Exponential-polynomial models}
\label{SEC exponential-polynomial models}

We now proceed to consider an example of a non-normalized parametric model, one of the major motivations to this work. In particular, let
\begin{align} \label{expon polyn model}
p_\vartheta(x) = C(\vartheta)^{-1} \exp\big( \vartheta_1 x + \dotso + \vartheta_d x^d \big), \quad x > 0,
\end{align}
$\vartheta = (\vartheta_1, \dots, \vartheta_d) \in \R^{d - 1} \times (- \infty, 0) = \Theta$, where
\begin{align*}
C(\vartheta) = \int_0^\infty \exp\big( \vartheta_1 x + \dotso + \vartheta_d x^d \big) \, \mathrm{d}x.
\end{align*}
These density functions correspond to a so-called exponential-polynomial model, which constitutes a special type of exponential family. It is trivial to see that these density functions obey the regularity assumptions (R1) - (R3), and also not hard to verify that the regularity conditions stated by \cite{BE:2018} (as summarized in Section \ref{SEC the new estimators}) are satisfied. Thus, we can first of all note, as a corollary to Theorem 3 of \cite{BE:2018}, the following characterization
\begin{corollary} \label{COR characterization expon polyn model}
	A positive random variable $X$ with $\E X^d < \infty$ follows the exponential-polynomial model in (\ref{expon polyn model}) if, and only if, the distribution function $F_X$ of $X$ satisfies
	\begin{align*}
	F_X(t)
	= \E \left[ - \left( \sum_{k = 1}^{d} k \vartheta_k X^{k - 1} \right) \min\{ X, t \} \right], \quad t > 0.
	\end{align*}
\end{corollary}
This is the characterization which underlies our new estimation method as constructed in Section \ref{SEC the new estimators}. \\

Notice that $C(\vartheta)$ cannot be written in a closed form, so maximum likelihood estimators are not readily available for the model in (\ref{expon polyn model}). Using the method of holonomic gradient descent, introduced by \cite{NNNOSTT:2011}, \cite{HT:2016} identify a differential equation which allows to numerically calculate $C(\vartheta)$ and its derivatives, and thus to get an approximation of the ML estimator. In our simulations, however, we focus on methods that do not try to approximate $C(\vartheta)$ numerically, but get rid of the normalization constant altogether. Namely, we consider our new method and compare it to the well-known score matching approach of \cite{H:2007}, in generalization of his method introduced in \cite{H:2005}, as well as to the noise-contrastive estimation technique of \cite{GH:2012}. In the case of non-negative, univariate observations, the score matching approach boils down to finding the minimum of
\begin{align*}
\widetilde{J}_{NN} (\vartheta)
&= \frac{1}{n} \sum_{j=1}^{n} \left[ 2 X_j \cdot \frac{p_\vartheta^\prime(X_j)}{p_\vartheta(X_j)} + \frac{p_\vartheta^{\prime\prime}(X_j)}{p_\vartheta(X_j)} \, X_j^2 - \frac{1}{2} \cdot \frac{\big(p_\vartheta^\prime(X_j)\big)^2}{\big(p_\vartheta(X_j)\big)^2} \, X_j^2 \right] \\
&= \frac{1}{n} \sum_{j=1}^{n} \left[ \left( \sum_{k = 1}^{d} k (k + 1) \vartheta_k X_j^k \right) + \frac{1}{2} \left( \sum_{k = 1}^{d} k \vartheta_k X_j^k \right)^2 \right]
\end{align*}
[see Section 3 of \cite{H:2007}], where $X_1, \dots, X_n$ are i.i.d. random variable with $X_1 \sim p_{\vartheta_0}$, for some unknown $\vartheta_0 \in \Theta$. Clearly, the quantity does not rely on $C(\vartheta)$. As for the estimator constructed in this paper, fixing $q = 2$ and the weight $w(t) = e^{- a t}$, where $a > 0$ is a tuning parameter, we may calculate $\psi_{n, 2}(\vartheta) = \big\lVert \eta_n( \, \cdot \, , \vartheta) \big\rVert_{L^2}$ (see Sections \ref{SEC the new estimators} and \ref{SEC Existence and measurability}) explicitly as
\begin{align*}
\psi_{n, 2}(\vartheta)
&= \frac{1}{n^2} \sum_{j = 1}^{n} \left\{ \frac{e^{- a X_j}}{a} \left[ 2 \sum_{\ell = 1}^d \ell \vartheta_\ell X_j^{\ell} + 1 - \left( \sum_{\ell = 1}^d \ell \vartheta_\ell X_j^{\ell - 1} \right)^2 \left( \frac{2 X_j}{a} + \frac{2}{a^2} \right) \right] + \frac{2}{a^3} \left( \sum_{\ell = 1}^d \ell \vartheta_\ell X_j^{\ell - 1} \right)^2 \right\} \\
&\quad - \frac{2}{n^2} \sum_{1 \leq j < k \leq n} \left\{ \left( \sum_{\ell = 1}^d \ell \vartheta_\ell X_{(j)}^{\ell} + 1 \right) \left[ - \frac{e^{- a X_{(k)}}}{a} + \frac{e^{- a X_{(k)}}}{a^2} \sum_{\ell = 1}^d \ell \vartheta_\ell X_{(k)}^{\ell - 1} \right] \right. \\
&\left. \qquad\qquad\qquad\qquad + \left( \sum_{\ell = 1}^d \ell \vartheta_\ell X_{(k)}^{\ell - 1} \right) \left[ \frac{e^{- a X_{(j)}}}{a^2} \left( \sum_{\ell = 1}^d \ell \vartheta_\ell X_{(j)}^{\ell} \right) - \frac{X_{(j)} \, e^{- a X_{(j)}}}{a} - \frac{e^{- a X_{(j)}}}{a^2} \right. \right. \\ 
&\left. \left. \qquad\qquad\qquad\qquad\qquad\qquad\qquad\qquad~ - \left( \sum_{\ell = 1}^d \ell \vartheta_\ell X_{(j)}^{\ell - 1} \right) \frac{2}{a^3} \big( 1 - e^{- a X_{(j)}} \big) \right] \right\} ,
\end{align*}
where $X_{(1)} < \dotso < X_{(n)}$ are the ordered values $X_1, \dots, X_n$. This formula is notably more complicated than the one resulting from the score matching approach, but in the two-parameter setting we now turn to, both estimators can be calculated explicitly. More precisely, to keep the presentation clear, we intent to focus on a two parameter case, but in order not to end up with a Gaussian-type model, we consider $d = 3$ and fix $\vartheta_2 = 0$, thus effectively considering the model
\begin{align*}
p_\vartheta(x) = C(\vartheta_1, \vartheta_3)^{-1} \exp\big( \vartheta_1 x + \vartheta_3 x^3 \big), \quad x > 0, \qquad \vartheta_1 \in \R, \, \vartheta_3 \in (- \infty, 0).
\end{align*}
In this case, each time by solving a quadratic equation in $\vartheta_1$ and $\vartheta_3$ (which is obtained by simplifying the above quantities $\widetilde{J}_{NN}$ and $\psi_{n, 2}$ further), we obtain the estimators explicitly. The score matching estimators for $\vartheta = (\vartheta_1, \vartheta_3)$ are given as
\begin{align*}
\widehat{\vartheta}_n^{SM}
= \left( - \frac{2 m_1}{m_2} - \frac{3 m_4}{m_2} \cdot \frac{4 m_2 m_3 - 2 m_1 m_4}{3 (m_4)^2 - 3 m_2 m_6}, ~ \frac{4 m_2 m_3 - 2 m_1 m_4}{3 (m_4)^2 - 3 m_2 m_6} \right)
\end{align*}
where $m_{k} = \sum_{j = 1}^{n} X_j^k$, and our new estimators are
\begin{align*}
\widehat{\vartheta}_{n, 2}^{(a)}
= \left( \frac{\overline{\Psi}_n^{(3)} \overline{\Psi}_n^{(5)} - 2 \overline{\Psi}_n^{(2)} \overline{\Psi}_n^{(4)}}{4 \overline{\Psi}_n^{(1)} \overline{\Psi}_n^{(2)} - \big( \overline{\Psi}_n^{(3)} \big)^2}, ~ \frac{\overline{\Psi}_n^{(3)} \overline{\Psi}_n^{(4)} - 2 \overline{\Psi}_n^{(1)} \overline{\Psi}_n^{(5)}}{4 \overline{\Psi}_n^{(1)} \overline{\Psi}_n^{(2)} - \big( \overline{\Psi}_n^{(3)} \big)^2} \right) ,
\end{align*}
where
\begin{align*}
\overline{\Psi}_n^{(1)} 
&= \frac{2}{a^3} + \frac{1}{n^2} \sum_{j = 1}^n e^{- a X_{(j)}} \left( - \frac{2 X_{(j)}}{a^2} - \frac{2}{a^3} \big( 2n - 2j + 1 \big) \right) + \frac{2}{n^2} \sum_{1 \leq j < k \leq n} - \frac{X_{(j)}}{a^2} \big( e^{- a X_{(k)}} + e^{- a X_{(j)}} \big) , \\
\overline{\Psi}_n^{(2)} 
&= \frac{1}{n^2} \sum_{j = 1}^n \left\{ - \frac{18 X_{(j)}^5}{a^2} \, e^{- a X_{(j)}} + \frac{18 X_{(j)}^4}{a^3} \big( 1 - e^{- a X_{(j)}} \big) \right\} \\
&\quad + \frac{2}{n^2} \sum_{1 \leq j < k \leq n} \left\{ - \frac{9 X_{(j)}^3 X_{(k)}^2}{a^2} \big( e^{- a X_{(k)}} + e^{- a X_{(j)}} \big) + \frac{18 X_{(j)}^2 X_{(k)}^2}{a^3} \big( 1 - e^{- a X_{(j)}} \big) \right\} , \\
\overline{\Psi}_n^{(3)} 
&= \frac{1}{n^2} \sum_{j = 1}^n \left\{ - \frac{12 X_{(j)}^3}{a^2} \, e^{- a X_{(j)}} + \frac{12 X_{(j)}^2}{a^3} \big( 1 - e^{- a X_{(j)}} \big) \right\} \\
&\quad + \frac{2}{n^2} \sum_{1 \leq j < k \leq n} \left\{ \big( e^{- a X_{(k)}} + e^{- a X_{(j)}} \big) \left( - \frac{3 X_{(j)} X_{(k)}^2}{a^2} - \frac{3 X_{(j)}^3}{a^2} \right) + \frac{6}{a^3} \big( 1 - e^{- a X_{(j)}} \big) \big( X_{(j)}^2 + X_{(k)}^2 \big) \right\} , \\
\overline{\Psi}_n^{(4)} 
&= \frac{1}{n^2} \sum_{j = 1}^n e^{- a X_{(j)}} \left( \frac{2 X_{(j)}}{a} + \frac{2 (n - 2j + 1)}{a^2} \right)
+ \frac{2}{n^2} \sum_{1 \leq j < k \leq n} \frac{X_{(j)}}{a} \big( e^{- a X_{(k)}} + e^{- a X_{(j)}} \big) ,
\end{align*}
and
\begin{align*}
\overline{\Psi}_n^{(5)} 
&= \frac{2}{n^2} \sum_{1 \leq j < k \leq n} \left\{ \frac{3 X_{(j)}^3}{a} \, e^{- a X_{(k)}} + \frac{3 X_{(j)} X_{(k)}^2}{a} \, e^{- a X_{(j)}} + \frac{3 X_{(k)}^2}{a^2} \big( e^{- a X_{(j)}} - e^{- a X_{(k)}} \big) \right\} \qquad\qquad\qquad\quad~ \\
&\quad + \frac{1}{n^2} \sum_{j = 1}^n \frac{6 X_{(j)}^3}{a} \, e^{- a X_{(j)}} .
\end{align*}
Moreover, we consider the noise-contrastive estimators in the refined version of \cite{GH:2012} that generalizes the initial results of \cite{GH:2010}. The idea is motivated by a binary classification problem and proceeds to consider the unknown normalization constant as an additional parameter to be estimated. The objective function is constructed in such a way that it ensures that the obtained estimator for the normalization constant truly provides (in numerical approximation) a normalized density without any further constraints on the optimization. Following \cite{GH:2012}, we implement this technique as follows. Given the sample $X_1, \dots, X_n$, choose the noise sample size $T_n = \nu \cdot n$ (for some fixed $\nu \in \N$, in our case $\nu = 10$) and sample from the noise distribution (in our case, the exponential distribution with rate parameter $\lambda_n = n / \sum_{j = 1}^n X_j$) to obtain values $Y_1, \dots, Y_{T_n}$. Then, minimize the objective function
\begin{align*}
J(\vartheta_1, \vartheta_3, c)
&= \frac{1}{n} \sum_{j = 1}^n \log\Big( 1 + \nu \, \lambda_n \, \exp\big( - (\lambda_n + \vartheta_1) X_j - \vartheta_3 X_j^3 - c \big) \Big) \\
&\quad - \frac{\nu}{T_n} \sum_{k = 1}^{T_n} \log\left( 1 - \frac{1}{1 + \nu \, \lambda_n \, \exp\big( - (\lambda_n + \vartheta_1) Y_k - \vartheta_3 Y_k^3 - c \big)} \right)
\end{align*}
to obtain an estimator $\widehat{\vartheta}_n^{NC}$ for the unknown parameters $(\vartheta_1^{(0)}, \vartheta_3^{(0)})$ as well as for the logarithm of the inverse of the normalization constant. In our simulations we used the 'L-BFGS-B'-method, which we have also applied in previous examples, for this optimization [with initial values $(0, -0.1, 0)$ and with the second parameter constrained to the negative numbers].

\begin{table}[ht!]
	\centering
	\setlength{\tabcolsep}{.4mm}
	\begin{tabular}{c||c||c|c||c|c|c|c|c|c}
		~~~$\vartheta_0$~~~ & ~~~$n$~~~ & ~~$\widehat{\vartheta}_n^{SM}$~~ & ~~$\widehat{\vartheta}_n^{NC}$~~ & ~~$\widehat{\vartheta}_{n, 2}^{(0.25)}$~~ & ~~$\widehat{\vartheta}_{n, 2}^{(0.5)}$~~ & ~~$\widehat{\vartheta}_{n, 2}^{(1)}$~~ & ~~$\widehat{\vartheta}_{n, 2}^{(2)}$~~ & ~~$\widehat{\vartheta}_{n, 2}^{(3)}$~~ & ~~$\widehat{\vartheta}_{n, 2}^{(5)}$~~ \\
		\hline
		& 10 & 1.4804 & \it (0.4925) & 0.4559 & 0.3977 & 0.3223 & \bf 0.2779 & 0.3353 & 0.7036 \\[1pt]
		\cline{3-10}
		& & -0.0619 & \it (-0.0268) & -0.0248 & -0.0223 & -0.0192 & \bf -0.0176 & -0.0206 & -0.0387 \\[1pt]
		\cline{2-10}
		& 25 & 0.5724 & \it (0.1589) & 0.1531 & 0.132 & 0.1067 & \bf 0.0903 & 0.096 & 0.1651 \\[1pt]
		\cline{3-10}
		& & -0.0221 & \it (-0.0084) & -0.008 & -0.0072 & -0.0062 & \bf -0.0057 & -0.0062 & -0.01 \\[1pt]
		\cline{2-10}
		${\, 1 \, \choose \, -0.05 \,}$ & 50 & 0.2982 & 0.074 & 0.0725 & 0.062 & 0.0495 & \bf 0.0406 & 0.0418 & 0.0641 \\[1pt]
		\cline{3-10}
		& & -0.0111 & \it (-0.0039) & -0.0038 & -0.0034 & -0.0029 & \bf -0.0026 & -0.0028 & -0.0042 \\[1pt]
		\cline{2-10}
		& 100 & 0.1521 & \it (0.0344) & 0.0339 & 0.029 & 0.0233 & \bf 0.0199 & 0.0215 & 0.0335 \\[1pt]
		\cline{3-10}
		& & -0.0056 & \it (-0.0019) & -0.0018 & -0.0016 & -0.0014 & \bf -0.0013 & -0.0014 & -0.0022 \\[1pt]
		\cline{2-10}
		& 200 & 0.079 & \it (0.0172) & 0.017 & 0.0145 & 0.0117 & \bf 0.0101 & 0.0109 & 0.0162 \\[1pt]
		\cline{3-10}
		& & -0.0029 & \it (-0.001) & -0.0009 & -0.0008 & \bf -0.0007 & \bf -0.0007 & \bf -0.0007 & -0.0011
		\\[1pt]
		\hline
		& 10 & 4.3787 & 0.9147 & 0.9046 & 0.8328 & 0.707 & 0.5174 & 0.3913 & \bf 0.253 \\[1pt]
		\cline{3-10}
		& & -1.9839 & -0.7007 & -0.6992 & -0.6667 & -0.6091 & -0.5203 & -0.4607 & \bf -0.4004 \\[1pt]
		\cline{2-10}
		& 25 & 1.6908 & 0.3153 & 0.3014 & 0.2722 & 0.2237 & 0.157 & 0.1166 & \bf 0.0738 \\[1pt]
		\cline{3-10}
		& & -0.6081 & -0.2038 & -0.1984 & -0.1868 & -0.1673 & -0.1407 & -0.1254 & \bf -0.1127 \\[1pt]
		\cline{2-10}
		${\, 0 \, \choose \, -0.5 \,}$ & 50 & 0.892 & 0.154 & 0.1456 & 0.131 & 0.1077 & 0.0771 & 0.0592 & \bf 0.0402 \\[1pt]
		\cline{3-10}
		& & -0.2963 & -0.0951 & -0.0917 & -0.0861 & -0.0771 & -0.0657 & -0.0595 & \bf -0.0549 \\[1pt]
		\cline{2-10}
		& 100 & 0.4785 & 0.0729 & 0.0691 & 0.0617 & 0.0499 & 0.0344 & 0.025 & \bf 0.0148 \\[1pt]
		\cline{3-10}
		& & -0.1512 & -0.0443 & -0.043 & -0.0402 & -0.0358 & -0.0301 & -0.0269 & \bf -0.0243 \\[1pt]
		\cline{2-10}
		& 200 & 0.2571 & 0.0375 & 0.0351 & 0.0314 & 0.0255 & 0.0177 & 0.0129 & \bf 0.0075 \\[1pt]
		\cline{3-10}
		& & -0.0789 & -0.0223 & -0.0213 & -0.0199 & -0.0177 & -0.0149 & -0.0133 & \bf -0.0119 \\[1pt]
		\hline
		& 10 & 8.5068 & 1.7608 & 1.7578 & 1.6817 & 1.5396 & 1.2927 & 1.0904 & \bf 0.7944 \\[1pt]
		\cline{3-10}
		& & -14.5606 & -5.1057 & -5.1765 & -5.0482 & -4.806 & -4.3777 & -4.0194 & \bf -3.4851 \\[1pt]
		\cline{2-10}
		& 25 & 3.2683 & 0.6169 & 0.5935 & 0.5616 & 0.5036 & 0.4084 & 0.3361 & \bf 0.2394 \\[1pt]
		\cline{3-10}
		& & -4.3203 & -1.481 & -1.4566 & -1.4099 & -1.3245 & -1.1834 & -1.0761 & \bf -0.9369 \\[1pt]
		\cline{2-10}
		${\, -0.5 \, \choose \, -3 \,}$ & 50 & 1.7089 & 0.2881 & 0.2739 & 0.2577 & 0.2288 & 0.1832 & 0.1499 & \bf 0.1066 \\[1pt]
		\cline{3-10}
		& & -2.0575 & -0.6602 & -0.6428 & -0.62 & -0.5793 & -0.5151 & -0.4688 & \bf -0.4119 \\[1pt]
		\cline{2-10}
		& 100 & 0.9238 & 0.1453 & 0.1376 & 0.1293 & 0.1148 & 0.0921 & 0.0758 & \bf 0.0545 \\[1pt]
		\cline{3-10}
		& & -1.0552 & -0.3193 & -0.3114 & -0.2999 & -0.2798 & -0.2486 & -0.2266 & \bf -0.1996 \\[1pt]
		\cline{2-10}
		& 200 & 0.4917 & 0.073 & 0.0674 & 0.0634 & 0.0564 & 0.0457 & 0.038 & \bf 0.0282 \\[1pt]
		\cline{3-10}
		& & -0.5441 & -0.1588 & -0.1518 & -0.1463 & -0.1367 & -0.1222 & -0.1122 & \bf -0.1002 \\[1pt]
	\end{tabular}
	\caption{Approximated biases calculated with 100,000 Monte Carlo samples.} \label{EXPON POLYN distribution Biases}
\end{table}

As in the previous simulations, we approximate bias and MSE of the competing estimators. The results are presented in Tables \ref{EXPON POLYN distribution Biases} and \ref{EXPON POLYN distribution MSE}. In the tables, for each underlying parameter $\vartheta_0 = \big( \vartheta_1^{(0)}, \vartheta_3^{(0)} \big)$ and each sample size, the first line corresponds to the bias/MSE of the $\vartheta_1$-parameter, while the second line corresponds to the $\vartheta_3$-parameter.

\begin{table}[ht!]
	\centering
	\setlength{\tabcolsep}{.4mm}
	\begin{tabular}{c||c||c|c||c|c|c|c|c|c}
		~~~$\vartheta_0$~~~ & ~~~$n$~~~ & ~~$\widehat{\vartheta}_n^{SM}$~~ & ~~$\widehat{\vartheta}_n^{NC}$~~ & ~~$\widehat{\vartheta}_{n, 2}^{(0.25)}$~~ & ~~$\widehat{\vartheta}_{n, 2}^{(0.5)}$~~ & ~~$\widehat{\vartheta}_{n, 2}^{(1)}$~~ & ~~$\widehat{\vartheta}_{n, 2}^{(2)}$~~ & ~~$\widehat{\vartheta}_{n, 2}^{(3)}$~~ & ~~$\widehat{\vartheta}_{n, 2}^{(5)}$~~ \\
		\hline
		& 10 & 5.9091 & \it (1.5668) & 1.4652 & 1.3332 & \bf 1.2118 & 1.5497 & 3.2052 & 25.4377 \\[1pt]
		\cline{3-10}
		& & 0.012 & \it (0.0041) & 0.0038 & 0.0035 & \bf 0.0031 & 0.0037 & 0.0068 & 0.0434 \\[1pt]
		\cline{2-10}
		& 25 & 1.1954 & \it (0.281) & 0.2982 & 0.2726 & \bf 0.2568 & 0.357 & 0.7326 & 4.1544 \\[1pt]
		\cline{3-10}
		& & 0.0019 & \textit{\textbf{(0.0006)}} & \bf 0.0006 & \bf 0.0006 & \bf 0.0006 & 0.0007 & 0.0015 & 0.0082 \\[1pt]
		\cline{2-10}
		${\, 1 \, \choose \, -0.05 \,}$ & 50 & 0.4758 & 0.1143 & 0.1228 & 0.1133 & \bf 0.109 & 0.1558 & 0.3253 & 1.6347 \\[1pt]
		\cline{3-10}
		& & 0.0007 & \textit{\textbf{(0.0002)}} & \bf 0.0002 & \bf 0.0002 & \bf 0.0002 & 0.0003 & 0.0006 & 0.0033 \\[1pt]
		\cline{2-10}
		& 100 & 0.2113 & \it (0.0514) & 0.055 & 0.0511 & \bf 0.05 & 0.0735 & 0.1541 & 0.7475 \\[1pt]
		\cline{3-10}
		& & 0.0003 & \textit{\textbf{(0.0001)}} & \bf 0.0001 & \bf 0.0001 & \bf 0.0001 & \bf 0.0001 & 0.0003 & 0.0015 \\[1pt]
		\cline{2-10}
		& 200 & 0.1012 & \it (0.0249) & 0.0262 & 0.0245 & \bf 0.0241 & 0.0358 & 0.0752 & 0.3573 \\[1pt]
		\cline{3-10}
		& & 0.0001 & \textit{\textbf{(0.0)}} & \bf 0.0 & \bf 0.0 & \bf 0.0 & 0.0001 & 0.0001 & 0.0007 \\[1pt]
		\hline
		& 10 & 42.1253 & 4.872 & 5.79 & 5.5517 & 5.188 & \bf 4.8596 & 5.0047 & 6.6628 \\[1pt]
		\cline{3-10}
		& & 15.2937 & 3.2687 & 3.4841 & 3.368 & 3.1697 & 2.899 & \bf 2.7841 & 3.012 \\[1pt]
		\cline{2-10}
		& 25 & 7.3114 & 1.1328 & 1.2154 & 1.1777 & \bf 1.1311 & 1.1409 & 1.275 & 1.9169 \\[1pt]
		\cline{3-10}
		& & 1.1392 & 0.2851 & 0.2938 & 0.2836 & 0.269 & \bf 0.2614 & 0.2814 & 0.41 \\[1pt]
		\cline{2-10}
		${\, 0 \, \choose \, -0.5 \,}$ & 50 & 2.6352 & 0.493 & 0.4993 & 0.4874 & \bf 0.4765 & 0.4997 & 0.5766 & 0.8948 \\[1pt]
		\cline{3-10}
		& & 0.3115 & 0.0934 & 0.0929 & 0.0905 & \bf 0.0878 & 0.0907 & 0.104 & 0.166 \\[1pt]
		\cline{2-10}
		& 100 & 1.1005 & 0.2248 & 0.2221 & 0.2175 & \bf 0.2144 & 0.2292 & 0.2688 & 0.4246 \\[1pt]
		\cline{3-10}
		& & 0.1097 & 0.0362 & 0.0358 & 0.035 & \bf 0.0345 & 0.037 & 0.0439 & 0.0736 \\[1pt]
		\cline{2-10}
		& 200 & 0.509 & 0.1099 & 0.1066 & 0.1046 & \bf 0.1038 & 0.112 & 0.132 & 0.2093 \\[1pt]
		\cline{3-10}
		& & 0.0456 & 0.0164 & 0.0159 & \bf 0.0156 & \bf 0.0156 & 0.017 & 0.0205 & 0.0351 \\[1pt]
		\hline
		& 10 & 158.341 & \bf 17.0727 & 21.1996 & 20.723 & 19.8873 & 18.6499 & 17.9361 & 17.9169 \\[1pt]
		\cline{3-10}
		& & 911.2872 & 171.7483 & 206.1279 & 202.4045 & 195.468 & 183.5396 & 174.0526 & \bf 161.9207 \\[1pt]
		\cline{2-10}
		& 25 & 26.1605 & 3.9728 & 4.2928 & 4.2173 & 4.0966 & \bf 3.9673 & 3.975 & 4.3622 \\[1pt]
		\cline{3-10}
		& & 56.2301 & 14.088 & 14.6718 & 14.3849 & 13.8935 & 13.2242 & \bf 12.952 & 13.5457 \\[1pt]
		\cline{2-10}
		${\, -0.5 \, \choose \, -3 \,}$ & 50 & 9.2065 & 1.6764 & 1.7049 & 1.6816 & 1.6482 & \bf 1.6295 & 1.6687 & 1.9025 \\[1pt]
		\cline{3-10}
		& & 14.5421 & 4.2357 & 4.2429 & 4.1757 & 4.0717 & \bf 3.9785 & 4.0302 & 4.5485 \\[1pt]
		\cline{2-10}
		& 100 & 3.7881 & 0.7814 & 0.7736 & 0.7647 & 0.7531 & \bf 0.7527 & 0.7795 & 0.9053 \\[1pt]
		\cline{3-10}
		& & 4.9664 & 1.6415 & 1.6254 & 1.6052 & 1.5777 & \bf 1.5715 & 1.6264 & 1.9134 \\[1pt]
		\cline{2-10}
		& 200 & 1.719 & 0.3773 & 0.3638 & 0.3601 & \bf 0.3559 & 0.3585 & 0.374 & 0.4388 \\[1pt]
		\cline{3-10}
		& & 1.9962 & 0.7272 & 0.7056 & 0.6988 & \bf 0.6914 & 0.6987 & 0.7334 & 0.883 \\[1pt]
	\end{tabular}
	\caption{Approximated MSE values calculated with 100,000 Monte Carlo samples.} \label{EXPON POLYN distribution MSE}
\end{table}

It is immediate that our new estimator and the noise-contrastive estimator outperform the score matching method distinctively over all tuning parameters, sample sizes, and parameter values for both the bias and the MSE, with the only exception being the parameter vector $(1, -0.05)$, for which the estimator $\widehat{\vartheta}_{n, 2}^{(5)}$ fares worse than the score matching approach in MSE terms. We propose as a very good compromise choice of the tuning parameter the use of $\widehat{\vartheta}_{n, 2}^{(1)}$ as an estimator. This particular estimator outperforms the score matching method by factors of (at least) $4$ in terms of MSE and also fares notably better in terms of the bias. It also outperforms the noise-contrastive estimation method uniformly, except for four instances in the MSE values (in three of which our method still performs better when another tuning parameter is chosen). The simulation in this non-normalized models conforms with the observation from previous examples that the new method fares remarkably well bias-wise. We also note that all of the estimators admit a large mean squared error for very small sample sizes, a behavior to be expected. From our simulations we conclude that the new estimation method is to be preferred clearly over the other approaches in this univariate setting of the exponential-polynomial models, but of course larger scale simulations involving different types of multi-parameter versions of the model would be needed to further strengthen this position [also, generalizations of the score matching technique, like \cite{YDS:2019}, could be taken into account]. One massive advantage of the score matching and noise-contrastive estimation approaches, however, is that they readily generalize to the multivariate situation, a generalization we were not (yet) able to establish for our approach (see the last paragraph of Section \ref{SEC notes and comments}).

\begin{remark} \label{RMK computational problems noise contrastive estimation}
	We observed in our simulations that the noise-contrastive estimators can run into computational problems when the exponentials in the objective function raise an overflow warning. A step by step analysis of the code suggests that for large noise sample sizes $T_n$ (that is, for large $\nu$) one tends to obtain some large values in the sample $Y_1, \dots, Y_{T_n}$ which are cubed in the exponential terms and thus become very (if not too) large. The behavior seems to appear more often for small parameter values $\vartheta_3^{(0)}$, but it seems to affect only single evaluations of the objective function during the optimization routine. We believe that most values for the noise-contrastive estimation approach in the table are intact and they also replicated when we reran the whole simulation, with a bit of an exception in the case of the parameter vector $(1, -0.05)$, where the values show a rather noticeable dependence on the initial value chosen for the optimization (though this does not happen for the other parameter values). Therefore, one possible ways to reduce the occurrence of overflows, which lies in choosing small initial values for the $\vartheta_3$-parameter in the optimization routine has an impact on the performance of the estimator. Another way out could be to adopt noise distributions with extremely short tails. It could prove useful to see if our observations replicate in other simulation studies. Note that no computational issues arise for the score matching and our new approach, where the estimators can be calculated explicitly.
\end{remark}

\section{Notes and comments}
\label{SEC notes and comments}

Note that there remain some problems for further research on our newly proposed estimators, the discussion or extension of which would be too extensive for this contribution. First, for all estimators we considered explicitly, we incorporate a tuning parameter '$a$' on which the performance depends strongly. It would be beneficial to have an adaptive choice of this parameter [see \cite{AS:2015}, and the refinements by \cite{T:2019}, who discuss such a method in the context of goodness-of-fit testing problems], probably adaptable to which criterion (minimal bias etc.) the estimator should satisfy. In the context of deriving results for $a \to \infty$, we obtained another consistent estimator for the Rayleigh parameter, and it would be interesting to see if such results can be derived for other distributions. Also, we have not used in practice the flexibility gained by providing all results for the general $L^q$-spaces, but restricted our attention to the case $q = 2$, mostly because of the explicit formulae obtainable in that case. If no closed formula for $\psi_{n, q}$ is feasible, either because of the use of some $q \neq 2$ or because some advanced weight function $w$ is chosen, the integral in $\psi_{n, q}$ has to be solved numerically which could lead to a computationally highly demanding procedure overall. As for the choice of a specific weight function $w$, to our best knowledge there exist no theoretical results which favor specific choices over others. Considering the vast amount of weighted $L^2$-statistics put to use in goodness-of-fit testing problems, it seems we cannot hope for general results in that direction. As such, the choice of the weight function provides some flexibility, but without clear guidance to satisfy specific objectives other than $\psi_{n, q}$ being calculable explicitly. \\

We have proven in a quite usual setting the consistency of our estimators. Surely, a limit theorem of the type
\begin{align*}
s(n) \big( \widehat{\vartheta}_{n, q} - \vartheta_0 \big) \stackrel{d}{\longrightarrow} \mathcal{P},
\end{align*}
where $s(n) \longrightarrow \infty$, as $n \to \infty$, and where $\mathcal{P}$ is some limit distribution (e.g. the normal distribution) is desirable. Such a result would pave the way for constructing confidence regions for the true parameter based on our method. The main hurdle in direct approaches of proving such a limit results, like some Taylor expansion or methods from empirical process theory, is that the terms involved in such calculations become too complicated and make the endeavor appear impractical to us. One hope is that, since \cite{BBDGM:2019} provide limit results for special classes of Stein discrepancy-based estimators, the interpretation of our estimation method in terms of the feature Stein discrepancy might at some point lead to advances. \\

Moreover, a larger-scale simulation study, involving more underlying parameters, sample sizes, and distributions could provide further insight into the estimation method. Improvements from a numerical point of view would, of course, benefit the approach. From a theoretical perspective, an important step in this last direction is to study whether the minimization method that is used in cases where the estimators cannot be calculated explicitly will always find a global minimum, or if not, in which situations it is likely to get stuck in some local minimum. \\

Note that \cite{BE:2018} also give characterization results for density functions on bounded intervals or on the whole real line. These can be used to construct similar estimation methods in the corresponding cases. To sketch the idea in the case of parametric models on the whole real line, assume that the support of each density function $p_{\vartheta}$ in $\mathfrak{P}_\Theta$ is the whole real line (and that some mild regularity conditions hold). Let $\widetilde{X}$ be a real-valued random variable with
\begin{align*}
\E \left[ \left| \frac{p_\vartheta^{\prime}(\widetilde{X})}{p_\vartheta(\widetilde{X})} \right| \Big( |\widetilde{X}| + 1 \Big) \right] < \infty, \quad \vartheta \in \Theta,
\end{align*}
and consider
\begin{align*}
\widetilde{\eta}(t, \vartheta) = \E\left[ \frac{p^\prime_\vartheta(\widetilde{X})}{p_\vartheta(\widetilde{X})} \, \big( t - \widetilde{X} \big) \, \mathds{1}\{ \widetilde{X} \leq t \} \right] - F_{\widetilde{X}}(t)
\end{align*}
for $(t, \vartheta) \in \R \times \Theta$. Then, similar to our elaborations in Section \ref{SEC the new estimators}, Theorem 4.1 of \cite{BE:2018} shows that $\widetilde{X} \sim p_{\vartheta_0}$ if, and only if, $\widetilde{\eta}(t, \vartheta_0) = 0$ for every $t \in \R$. Therefore, if, initially, $\widetilde{X} \sim p_{\vartheta_0}$, then $\lVert \widetilde{\eta}(\cdot \, , \vartheta) \rVert_{L^q} = 0$ if, and only if, $\vartheta = \vartheta_0$. Here, $L^q = L^q\big( \R, \mathcal{B}^1, \widetilde{w}(t) \, \mathrm{d}t \big)$, $1 \leq q < \infty$, with a positive weight function $\widetilde{w}$ satisfying
\begin{align*}
\int_\R \big( |t|^q + 1 \big) \, \widetilde{w}(t) \, \mathrm{d}t < \infty.
\end{align*}
Thus, with
\begin{align*}
\widetilde{\eta}_n(t, \vartheta) = \frac{1}{n} \sum_{j=1}^{n} \frac{p^\prime_\vartheta(\widetilde{X}_j)}{p_\vartheta(\widetilde{X}_j)} \, \big( t - \widetilde{X}_j \big) \, \mathds{1}\{ \widetilde{X}_j \leq t \} - \frac{1}{n} \sum_{j=1}^{n} \mathds{1}\{ \widetilde{X}_j \leq t \} ,
\end{align*}
a reasonable estimator for $\vartheta_0$ is
\begin{align*}
\widetilde{\vartheta}_{n, q}
= \arg\min\big\{ \lVert \widetilde{\eta}_n(\cdot \, , \vartheta) \rVert_{L^q} \, | \, \vartheta \in \Theta \big\}.
\end{align*}
Apparently, once we switch to density function supported by the whole real line, the characterization result due to \cite{BE:2018}, and thus our estimator, have slightly different forms, but using the results from Section \ref{SEC Existence and measurability}, we could still prove existence and measurability for this type of estimator, and give a formal definition as in (\ref{Estimator, formal defintion}). Moreover, a classical proof via the law of large numbers for random elements in separable Banach spaces and the Arzel\`{a}-Ascoli theorem [considering the modulus of continuity, as employed by \cite{B:1968}] yields the convergence results from Lemma \ref{LEMMA convergence of psi_n,q} for $\widetilde{\psi}_{n, q} = \lVert \widetilde{\eta}_n(\cdot \, , \vartheta) \rVert_{L^q}$, but with all convergences only in probability. That result can then be used to derive consistency as in Theorem \ref{THM consistency}, again with all convergences only in probability. However, choosing a fixed (i.e. parameter-independent) weight function on $\R$ with a mere scale-tuning, as we employ it throughout (using the weight $t \mapsto e^{- a t}$), appears not to be sufficient to account for the possible location-dependence of the model. Thus, in simulations (for instance with the Cauchy distribution) the problem, to us, seems empirically more involved and is therefore not addressed in the work at hand. \\

Still, we deem it possible to apply our new type of estimator to models which are supported by any connected subset of $\R$ as indicated in the previous lines. Of course, the next question which forces itself on us is whether a similar method can be devised for multivariate models. Here the frontiers are somewhat blurry: The Stein density approach identity which appears at the beginning of Section \ref{SEC the new estimators} is not yet fully understood in the multivariate case [as stated in Remark 1.1 by \cite{LRS:2017:2}], and the characterizations derived by \cite{BE:2018} rely on further calculations, the generalization of which is not immediate. Thus, we have to state at this point that, to us, it is an open question how a generalization to the multivariate setting could look like (with no clear indication of it being possible at all).

\appendix

\section{Additional material for Section \ref{SEC Existence and measurability}} 
\label{Appendix existence section}

\begin{remark}{[Comments on Theorem \ref{THM measurable selection}]} \label{RMK comments to selection theorem}
	There is another result which gives measurable selections without the completeness assumption on the probability space [as provided by \cite{BP:1973}], but it requires $\sigma$-compactness of the parameter space, thus essentially reducing the study to euclidean parameters (a Banach space is $\sigma$-compact if, and only if, it is of finite dimension, which follows easily from Baire's category theorem). Of course this is enough for our purposes, but currently the interest in statistical inference for infinite dimensional models grows remarkably. Hence if a statistician was to investigate measurability of an estimator for some infinite dimensional quantity, she would have to resort to a result in the generality of Theorem \ref{THM measurable selection}. Another reason for us to build on Theorem \ref{THM measurable selection} is that other measurability results known to us do not quite fit the construction of our estimators. For instance, \cite{S:1970} considers minimum discrepancy estimators, where discrepancies are (certain) functions on the Cartesian product of a suitable set of probability measures with itself. It is (formally) not possible to identify such a set of probability measures in our setting, as we ought to introduce the empirical distribution of a sample into the discrepancy function, while only considering parametric distributions with a continuously differentiable density. Even though we believe this to be a purely formal issue which might be resolved to render results from \cite{S:1970} applicable, additional caution is needed that Theorem \ref{THM measurable selection} does not require. Likewise, the setting considered by \cite{P:1969} does not cover our estimators. \\
	
	Note that since completing (the $\sigma$-field of) an underlying probability space does not interfere with measurability properties of random maps, nor does it meddle with push-forward measures, the corresponding assumption in Theorem \ref{THM measurable selection} is no restriction. If $\mathfrak{S}$ is a complete, separable metric space and the map $\Gamma$ from Theorem \ref{THM measurable selection} takes compact subsets of $\mathfrak{S}$ as values, the condition imposed on the graph is equivalent to $\Gamma$ being measurable with respect to the Borel-$\sigma$-field generated by the Hausdorff topology [see Theorems III.2 and III.30 by \cite{CV:1977}]. Likewise, if $\mathfrak{S}$ is a locally compact, complete, separable metric space and $\Gamma$ maps into the closed subsets of $\mathfrak{S}$, the condition is equivalent to $\Gamma$ being measurable with respect to the Borel-$\sigma$-field generated by the Fell topology [this can be proven using results from \cite{B:1993} and \cite{CV:1977}]. \\
\end{remark}
\noindent
\textbf{\textit{Proof of Lemma \ref{LEMMA measurability eta_n}}}.
First recall the following lemma on product-measurability, the proof of which is an easy exercise.
\begin{lemma} \label{LEMMA product measurability}
	Let $(S, \mathcal{A}, \mu)$ be a measure space, $I \subset \R$ an open interval, and let $(\mathcal{T}, \mathcal{O}_\mathcal{T})$ be a topological vector space. Furthermore, let $h : S \times I \to \mathcal{T}$ be a map such that
	\begin{itemize}
		\item $s \mapsto h(s, x)$ is $\big( \mathcal{A}, \mathcal{B}(\mathcal{T}) \big)$-measurable for every $x \in I$, and
		\item $x \mapsto h(s, x)$ is right-continuous for every $s \in S$.
	\end{itemize}
	Then $h$ is $\big( \mathcal{A} \otimes \mathcal{B}(I), \mathcal{B}(\mathcal{T}) \big)$-measurable.
\end{lemma}
Notice that for any fixed $(t, \vartheta) \in (0, \infty) \times \Theta$ the map $\omega \mapsto \eta_n(\omega, t, \vartheta)$ is $(\mathcal{F}, \mathcal{B}^1)$-measurable, and for any fixed $(\omega, t) \in \Omega \times (0, \infty)$ the map $\vartheta \mapsto \eta_n(\omega, t, \vartheta)$ is continuous. By a statement analogous to Lemma \ref{LEMMA product measurability} [see for instance Lemma III.14 by \cite{CV:1977}], $(\omega, \vartheta) \mapsto \eta_n(\omega, t, \vartheta)$ is $\big( \mathcal{F} \otimes \mathcal{B}(\Theta), \mathcal{B}^1 \big)$-measurable for fixed $t > 0$. Since $t \mapsto \eta_n(\omega, t, \vartheta)$ is continuous for fixed $(\omega, \vartheta) \in \Omega \times \Theta$, Lemma \ref{LEMMA product measurability} implies that $\eta_n$ is $\big( \mathcal{F} \otimes \mathcal{B}(0, \infty) \otimes \mathcal{B}(\Theta), \mathcal{B}^1 \big)$-measurable. Consequently, the maps
\begin{align*}
(\omega, \vartheta) \mapsto \big\langle \eta_n(\omega, \cdot \, , \vartheta), g \big\rangle_{L^q}
\end{align*}
are measurable for every $g \in L^{q^\prime}$ by Fubini's theorem, and since $L^q$ is a separable Banach space, the mapping $(\omega, \vartheta) \mapsto \eta_n(\omega, \cdot\, , \vartheta)$ is $\big( \mathcal{F} \otimes \mathcal{B}(\Theta), \mathcal{B}(L^q) \big)$-measurable [cf. Corollary 1.1.2 of \cite{HNVW:2016}]. \hfill\ensuremath{\square}

\begin{remark}{[$\Gamma_{n, q}$ from (\ref{set of approximate estimators}) is closed]} \label{RMK closedness of Gamma_n,q}
	Note that (R1) and Fatou's lemma imply the lower semi-continuity of the map $\vartheta \mapsto \psi_{n, q}(\omega, \vartheta)$. Thus if $\vartheta^{(k)} \in \Gamma_{n, q}(\omega)$, $k \in \N$, converges (with respect to the metric in $\Theta$) to $\vartheta^* \in \Theta$ as $k \to \infty$, then
	\begin{align*}
	\psi_{n, q}(\omega, \vartheta^*) 
	= \big\lVert \eta_n(\omega, \cdot \, , \vartheta^*) \big\rVert_{L^q}
	\leq \liminf_{k \, \to \, \infty} \big\lVert \eta_n\big(\omega, \cdot \, , \vartheta^{(k)}\big) \big\rVert_{L^q}
	\leq m_{n, q}(\omega) + \varepsilon_n(\omega) ,
	\end{align*}
	that is, $\vartheta^* \in \Gamma_{n, q}(\omega)$, so $\Gamma_{n, q}(\omega)$ is closed in $\Theta$ for every $\omega \in \Omega$. Hence we can note that if $\Theta$ is closed, and therefore locally compact [cf. p.42 of \cite{K:1968}] and complete, $\Gamma_{n, q}$ is a random element in the space of all closed subsets of $\Theta$ endowed with the Fell topology (see also Remark \ref{RMK comments to selection theorem}).
\end{remark}

\section{Additional material for Section \ref{SEC Consistency}}
\label{Appendix consistency section}

\noindent
\textbf{\textit{Proof of Lemma \ref{LEMMA convergence of psi_n,q}}}.
First note that for any non-empty closed subset $F$ of $K$,
\begin{align*}
\Big| \inf_{\vartheta \, \in \, F} \psi_{n, q}(\vartheta) - \inf_{\vartheta \, \in \, F} \psi_q(\vartheta) \Big|
\leq \sup_{\vartheta \, \in \, K} \big| \psi_{n, q}(\vartheta) - \psi_q(\vartheta) \big| ,
\end{align*}
so the second claim of Lemma \ref{LEMMA convergence of psi_n,q} follows from the first. For the first claim, let $K \neq \emptyset$ be a compact subset of $\Theta$. Note that
\begin{align} \label{estimate empirical process}
\sup_{\vartheta \, \in \, K} \big| \psi_{n, q}(\vartheta) - \psi_q(\vartheta) \big|
&\leq \sup_{\vartheta \, \in \, K} \big\lVert \eta_n(\, \cdot \, , \vartheta) - \eta(\, \cdot \, , \vartheta) \big\rVert_{L^q} \nonumber \\
&\leq C \cdot \sup_{\substack{\vartheta \, \in \, K \\ t \, > \, 0}} \left| \frac{1}{n} \sum_{j=1}^{n} \frac{p^\prime_{\vartheta}(X_j)}{p_\vartheta(X_j)} \, \min\{ X_j , t \} - \E\left[ \frac{p^\prime_{\vartheta}(X)}{p_\vartheta(X)} \, \min\{ X, t \} \right] \right| \nonumber \\
&~~~~+ C \cdot \sup_{t \, > \, 0} \left| \frac{1}{n} \sum_{j=1}^{n} \mathds{1}\{ X_j \leq t \} - F_{X}(t) \right| ,
\end{align}
where $C = \big( \int_0^\infty w(t) \, \mathrm{d}t \big)^{1/q}$. The second term on the right-hand side of (\ref{estimate empirical process}) converges to $0$ almost surely by the classical Glivenko-Cantelli theorem. For a function $f : (0, \infty) \to \R$ we write $\mathbb{P}_n f = \tfrac{1}{n} \sum_{j=1}^{n} f(X_j)$ and $\mathbb{P}^{X} f = \E \big[ f(X) \big]$. Then the first term on the right-hand side of (\ref{estimate empirical process}) can be written as
\begin{align} \label{glivenko-cantelli class}
\sup_{\substack{\vartheta \, \in \, K \\ t \, > \, 0}} \left| \frac{1}{n} \sum_{j=1}^{n} \frac{p^\prime_{\vartheta}(X_j)}{p_\vartheta(X_j)} \, \min\{ X_j , t \} - \E\left[ \frac{p^\prime_{\vartheta}(X)}{p_\vartheta(X)} \, \min\{ X, t \} \right] \right|
= \sup_{\substack{\vartheta \, \in \, K \\ t \, > \, 0}} \Big| \mathbb{P}_n f_{t, \vartheta} - \mathbb{P}^{X} f_{t, \vartheta} \Big|
= \sup_{f \, \in \, \mathcal{H}_{\Theta}} \Big| \mathbb{P}_n f - \mathbb{P}^{X} f \Big| ,
\end{align}
where $f_{t, \vartheta}(x) = \tfrac{p^\prime_\vartheta(x)}{p_\vartheta(x)} \cdot \min\{ x, t \}$, $x > 0$, is a measurable function for every $\vartheta \in K$ and $t > 0$, and where $\mathcal{H}_{\Theta} = \big\{ f_{t, \vartheta} \, \big| \, \vartheta \in K, \, t > 0 \big\}$ denotes the collection of all such functions. Note that the supremum in (\ref{glivenko-cantelli class}) is finite ($\mathbb{P}$-a.s.) by (R1), (\ref{Integrability condition for X in the characterization}), and (R3), and that the terms in (\ref{glivenko-cantelli class}) constitute measurable maps from $(\Omega, \mathcal{F})$ to $(\R, \mathcal{B}^1)$ by Theorem \ref{THM measurable selection}.

As is commonly done, we denote, for given functions $l, u : (0, \infty) \to \R$, by $[l, u]$ the set of all functions $f$ such that $l \leq f \leq u$ pointwise. An $\varepsilon$-bracket with respect to $L^1(\mathbb{P}^{X}) = L^1\big( (0, \infty),\, \mathcal{B}(0, \infty),\, \mathbb{P}^{X} \big)$ is one such set $[l, u]$ with $\lVert u - l \rVert_{L^1(\mathbb{P}^{X})} < \varepsilon$. The bracketing number $\mathcal{N}_{[\,]} \big( \varepsilon, \mathcal{H}_{\Theta}, L^1(\mathbb{P}^{X}) \big)$ of $\mathcal{H}_\Theta$ is the minimum number of $\varepsilon$-brackets needed to cover $\mathcal{H}_\Theta$. If the bracketing number of $\mathcal{H}_\Theta$ is finite for every $\varepsilon > 0$, then $\mathcal{H}_\Theta$ is a Glivenko-Cantelli class, that is, $\sup_{f \, \in \, \mathcal{H}_{\Theta}} \big| \mathbb{P}_n f - \mathbb{P}^{X} f \big| \longrightarrow 0$ almost surely [see Theorem 2.4.1 by \cite{VW:2000}], which, combined with (\ref{estimate empirical process}) and (\ref{glivenko-cantelli class}), implies the claim. Note that the result by \cite{VW:2000} is formulated to give convergence outer almost surely, but as we work on a complete probability space, the transition to an outer probability measure is not necessary (since we can provide enough measurability on a complete probability space and the notions of almost sure convergence and outer almost sure convergence agree).

Thus, to prove Lemma \ref{LEMMA convergence of psi_n,q}, it remains to show that the bracketing numbers of $\mathcal{H}_\Theta$ are finite. The following argument combines ideas from the classical Glivenko-Cantelli theorem and from Example 19.7 of \cite{V:1998}. Let $\varepsilon > 0$ be arbitrary, and set $\delta = \varepsilon^{1 / \alpha} \, (4 \, \E[ H(X) \, X ])^{-1 / \alpha}$, where $H$ and $\alpha$ are as in (R3). Since $K$ is compact there exist $\vartheta_1, \dots, \vartheta_m \in K$, $m = m_\varepsilon \in \N$, such that $\bigcup_{i=1}^{m} B_\delta(\vartheta_i) \supset K$. Additionally, since for each $i = 1, \dots, m$ the function
\begin{align*}
[0, \infty) \ni t \mapsto E_i(t)
= \E \left[ \left| \frac{p_{\vartheta_i}^\prime(X)}{p_{\vartheta_i}(X)} \right| \, \min\{ X, t \} \right]
\end{align*}
is continuous and monotonically increasing, and since it satisfies the limit relation $E_i(0) = \lim_{t \, \searrow \, 0} E_i(t) = 0$ as well as $E_i(\infty) = \lim_{t \, \nearrow \, \infty} E_i(t) = \E \left| \tfrac{p_{\vartheta_i}^\prime(X)}{p_{\vartheta_i}(X)} \, X \right| < \infty$, there exist $0 = t_0 < t_1 < \dotso < t_\ell = \infty$, $\ell = \ell_\varepsilon \in \N$, such that
\begin{align*}
E_i(t_j) - E_i(t_{j - 1}) < \varepsilon / 4
\end{align*}
for $j = 1, \dots, \ell$ and $i = 1, \dots, m$. Upon setting $f_{0, \vartheta}(x) = 0$, $f_{\infty, \vartheta}(x) = \tfrac{p_\vartheta^\prime(x)}{p_\vartheta(x)} \cdot x$, for $x > 0$ and $\vartheta \in K$, we define the brackets
\begin{align*}
\mathcal{H}_{i, j}
= \Big[ f_{t_{j - 1}, \vartheta_i} - \big| f_{t_{j}, \vartheta_i} - f_{t_{j - 1}, \vartheta_i} \big| - \delta^\alpha \cdot H^*, \, f_{t_{j - 1}, \vartheta_i} + \big| f_{t_{j}, \vartheta_i} - f_{t_{j - 1}, \vartheta_i} \big| + \delta^\alpha \cdot H^* \Big] ,
\end{align*}
for $j = 1, \dots, \ell$ and $i = 1, \dots, m$, where $H^*(x) = H(x) \cdot x$, $x > 0$. These brackets cover $\mathcal{H}_\Theta$. Indeed, if $\vartheta \in K$ and $t > 0$ are arbitrary, there exist $i \in \{1, \dots, m \}$ and $j \in \{ 1, \dots, \ell \}$ such that $\vartheta \in B_\delta(\vartheta_i)$ and $t_{j - 1} \leq t < t_j$, so $f_{t, \vartheta} \in \mathcal{H}_{i, j}$ since for every $x > 0$
\begin{align*}
\big| f_{t, \vartheta}(x) - f_{t_{j - 1}, \vartheta_i}(x) \big|
&\leq \big| f_{t, \vartheta}(x) - f_{t, \vartheta_i}(x) \big| + \big| f_{t, \vartheta_i}(x) - f_{t_{j - 1}, \vartheta_i}(x) \big| \\
&= \left| \frac{p_\vartheta^\prime(x)}{p_\vartheta(x)} - \frac{p_{\vartheta_i}^\prime(x)}{p_{\vartheta_i}(x)} \right| \, \min\{ x, t \} + \left| \frac{p_{\vartheta_i}^\prime(x)}{p_{\vartheta_i}(x)} \right| \Big( \min\{ x, t \} - \min\{ x, t_{j - 1} \} \Big) \\
&\leq H(x) \cdot x \cdot \big| \vartheta - \vartheta_i \big|^\alpha + \left| \frac{p_{\vartheta_i}^\prime(x)}{p_{\vartheta_i}(x)} \right| \Big( \min\{ x, t_j \} - \min\{ x, t_{j - 1} \} \Big) \\
&\leq \delta^\alpha \cdot H^*(x) + \big| f_{t_j, \vartheta_i}(x) - f_{t_{j - 1}, \vartheta_i}(x) \big| .
\end{align*}
Moreover, the brackets $\mathcal{H}_{i, j}$ are $\varepsilon$-brackets with respect to $L^1(\mathbb{P}^{X})$, as
\begin{align*}
\Big\lVert 2 \Big( \big| f_{t_j, \vartheta_i} - f_{t_{j - 1}, \vartheta_i} \big| + \delta^\alpha \cdot H^* \Big) \Big\rVert_{L^1(\mathbb{P}^{X})}
&= 2 \, \E \left[ \big| f_{t_j, \vartheta_i}(X) - f_{t_{j - 1}, \vartheta_i}(X) \big| + \delta^\alpha \cdot H(X) \, X \right] \\
&= 2 \big( E_i(t_j) - E_i(t_{j - 1}) \big) + \frac{\varepsilon}{2} \\
&< \varepsilon .
\end{align*}
Hence $\mathcal{N}_{[\,]} \big( \varepsilon, \mathcal{H}_{\Theta}, L^1(\mathbb{P}^{X}) \big) \leq m_\varepsilon \cdot \ell_\varepsilon < \infty$. \hfill\ensuremath{\square}

\vspace{5mm}
\noindent
\textbf{\textit{Proof of Remark \ref{RMK consistency in convex case}}}.
From Lemma \ref{LEMMA convergence of psi_n,q} we know that $\psi_{n, q}(\vartheta) \to \psi_q(\vartheta)$ $\mathbb{P}$-a.s., as $n \to \infty$, for each $\vartheta \in \Theta$. Since $\vartheta_0 \in \Theta^\circ$, there exists a $\delta > 0$ such that $B_{2 \delta}(\vartheta_0) \subset \Theta$. Then the closed ball $\overline{B} = \overline{B_{\delta}(\vartheta_0)}$ also lies in $\Theta$. Denote by $R = \partial B_{\delta}(\vartheta_0)$ the boundary of that ball. It follows from Lemma \ref{LEMMA convergence of psi_n,q} that
\begin{align*}
\inf_{\vartheta \, \in \, \overline{B}} \psi_{n, q}(\vartheta)
\longrightarrow
\inf_{\vartheta \, \in \, \overline{B}} \psi_q(\vartheta) = 0 \quad \text{and} \quad
\inf_{\vartheta \, \in \, R} \psi_{n, q}(\vartheta)
\longrightarrow
\inf_{\vartheta \, \in \, R} \psi_q(\vartheta) > 0,
\end{align*}
both $\mathbb{P}$-a.s., as $n \to \infty$, where the positiveness of the last term follows from (\ref{minimum is well separated}). Now, let $\varepsilon > 0$ and choose $n_0 = n_0(\varepsilon) \in \N$ such that
\begin{align*}
\mathbb{P}\Big( \inf_{\vartheta \, \in \, \overline{B}} \psi_{n, q}(\vartheta) + \varepsilon_n < \inf_{\vartheta \, \in \, R} \psi_{n, q}(\vartheta) \Big)
\geq 1 - \frac{\varepsilon}{2}, \quad n \geq n_0.
\end{align*}
Next, note that if $\inf_{\vartheta \, \in \, \overline{B}} \psi_{n, q}(\vartheta) + \varepsilon_n < \inf_{\vartheta \, \in \, R} \psi_{n, q}(\vartheta)$ then $\psi_{n, q}$ has a local minimum in $B_{\delta}(\vartheta_0)$ (since $\varepsilon_n > 0$) which, by strict convexity, is the unique global minimum. Additionally, we have 
\begin{align*}
	\inf_{\vartheta \, \in \, \overline{B}} \psi_{n, q}(\vartheta) + \varepsilon_n < \inf_{\vartheta \, \in \, \Theta \setminus B_\delta(\vartheta_0)} \psi_{n, q}(\vartheta).
\end{align*}
On the other hand, if we have the relation $\widehat{\vartheta}_{n, q} \in \Theta \setminus B_\delta(\vartheta_0)$, then
\begin{align*}
\inf_{\vartheta \, \in \, \Theta \setminus B_\delta(\vartheta_0)} \psi_{n, q}(\vartheta)
\leq \psi_{n, q}\big( \widehat{\vartheta}_{n, q} \big)
\leq \inf_{\vartheta \, \in \, \Theta} \psi_{n, q}(\vartheta) + \varepsilon_n
= \inf_{\vartheta \, \in \, \overline{B}} \psi_{n, q}(\vartheta) + \varepsilon_n.
\end{align*}
Consequently, for all $n \geq n_0$,
\begin{align*}
1 - \frac{\varepsilon}{2}
&\leq \mathbb{P}\Big( \inf_{\vartheta \, \in \, \overline{B}} \psi_{n, q}(\vartheta) + \varepsilon_n < \inf_{\vartheta \, \in \, R} \psi_{n, q}(\vartheta) \Big) \\
&\leq \mathbb{P}\Big( \inf_{\vartheta \, \in \, \overline{B}} \psi_{n, q}(\vartheta) + \varepsilon_n < \inf_{\vartheta \, \in \, \Theta \setminus B_\delta(\vartheta_0)} \psi_{n, q}(\vartheta) \Big) \\
&\leq \mathbb{P}\Big( \widehat{\vartheta}_{n, q} \in \overline{B} \Big).
\end{align*}
Since $\big\{ \mathbb{P}^{\widehat{\vartheta}_{n, q}} \, \big| \, n \leq n_0 \big\}$ is a finite set of measures, there exists a compact set $K \subset \R^d$ such that $\mathbb{P}\big( \widehat{\vartheta}_{n, q} \in K \big) \geq 1 - \tfrac{\varepsilon}{2}$ for all $n \leq n_0$. The set $K \cap \overline{B} \subset \Theta$ is a compact subset of $\R^d$ and thus also of $\Theta$, for a compact metric space is a compact subset of every metric space it embeds into continuously [see p.21, Theorem 3, of \cite{K:1968}]. By choice of the sets,
\begin{align*}
\mathbb{P}\Big( \widehat{\vartheta}_{n, q} \in K \cap \overline{B} \Big)
\geq 1 - \varepsilon,
\end{align*}
which is the claim. \hfill\ensuremath{\square}
\vspace{7mm}

\bibliography{references_parameter_estimation}   
\bibliographystyle{apalike}

\end{document}